\documentclass[11pt]{amsart}
\usepackage{amssymb,amsmath}
\usepackage{graphicx}
\usepackage[colorlinks=true]{hyperref}
\hypersetup{urlcolor=[rgb]{0.00,0.38,0.00},
linkcolor=[rgb]{0.00,0.00,0.70},
citecolor=[rgb]{0.85,0.00,0.00}
,colorlinks=[rgb]{0.03,1.00,0.00}
}
\long\def\symbolfootnote[#1]#2{\begingroup
\def\thefootnote{\fnsymbol{footnote}}\footnote[#1]{#2}\endgroup}
\font\bb=msbm10 at 12pt
\numberwithin{equation}{section}
\def\cqfd{\nopagebreak\hfill \vrule width 4pt height 4pt depth 1pt}
\font\bba=msbm10 \font\bbb=msbm8
\newfam\bbfam \textfont\bbfam=\bba \scriptfont\bbfam=\bbb
\def\bb{\fam\bbfam\bba} \def\N{{\bb N}} 
\def\u{{\mathcal U}}
\def\l{{\mathcal L}}

\def\H{{\mathcal H}}

\def\b{{\mathcal B}}

\def\Z{{\bb Z}} 
\def\R{{\bb R}}
\def\C{{\bb C}}

\def\re{{\mathcal Re}\ }
\def\im{{\mathcal Im}\ }

\def\1{1\!\!1}

\def\ds{\displaystyle}

\def\lam{\lambda}
\def\De{\Delta}
\def\dist{{\rm dist}}
\newtheorem{Th}{Theorem}[section] \newtheorem{Lemme}[Th]{Lemma}
\newtheorem{Cor}[Th]{Corollary} \newtheorem{Prop}[Th]{Proposition}

\newtheorem{Lemma}[Th]{Lemma}
\begin{document}
\title{Variations of Hausdorff Dimension in the exponential family}
\maketitle

\author{Guillaume Havard\symbolfootnote[2]{Laboratoire de Math\'ematiques (UMR 6620), universit\'e Blaise Pascal, campus universitaire des C\'ezeaux
63177 Aubi\`ere cedex, France (\tt{\href{mailto:guillaume.havard@math.univ-bpclemront.fr}{guillaume.havard@math.univ-bpclemront.fr}})},
\, Mariusz Urba\'nski\symbolfootnote[3]{Department of Mathematics, University of North Texas, P.O. 311430, Denton TX
76203-1430, USA,
(\tt{\href{mailto:urbanski@unt.edu}{urbanski@unt.edu},\, \href{http://www.math.unt.edu/~urbanski/}{http://www.math.unt.edu/$\sim$urbanski}}).
The research of Mariusz Urbanski is partially supported by the NSF grant DMS 0700831.
}\,
and\, Michel Zinsmeister\symbolfootnote[4]{Laboratoire de Math\'ematiques et Applications, Physique Math\'ematique d'Orl\'eans
(MAPMO - UMR6628) F\'ed\'eration Denis Poisson (FDP - FR2964) CNRS/Universit\'e d'Orl\'eans, B.P. 6759, 45067
Orl\'eans cedex 2, France (\tt{\href{mailto:Michel.Zinsmeister@math.cnrs.fr}{Michel.Zinsmeister@math.cnrs.fr},\,
\href{http://www.univ-orleans.fr/mapmo/membres/zins/}{http://www.univ-orleans.fr/mapmo/membres/zins/}})}}

\begin{abstract}
In this paper we deal with the following family  of exponential
maps $(f_\lambda:z\mapsto \lambda(e^z-1))_{\lambda\in
  [1,+\infty)}$. Denoting
$d(\lambda)$ the hyperbolic dimension of $f_\lambda$. It is proved in \cite{uz2} that the
function $\lambda\mapsto d(\lambda)$ is real analytic in $(1,+\infty)$, and in \cite{uz1}
that it is continuous in $[1,+\infty)$. In this paper we prove that this map is C$^1$
on $[1,+\infty)$, with $d'(1^+)=0$. Moreover we prove that depending on the value of
$d(1)$
$$
\left\{\begin{array} {rcll} d'(1+\varepsilon)
& \sim  & -\varepsilon^{2d(1)-2} & \mbox{if $d(1)
<\frac{3}{2}$,}\vspace{0.2cm}\\
|d'(1+\varepsilon)|& \lesssim & -\varepsilon\log\varepsilon & \mbox{if
  $d(1)
=\frac{3}{2},$}\vspace{0.2cm}\\
|d'(1+\varepsilon)|
& \lesssim & \varepsilon & \mbox{if $d(1)>\frac{3}{2}.$}
\end{array}\right.
$$
In particular, if $d(1)<\frac{3}{2}$, then there exists $\lambda_0>1$ such that $d(\lambda)<d(1)$
for any $\lambda\in(1,\lambda_0)$.
\end{abstract}


\tableofcontents

\section{Introduction}
\subsection{An overview of the problem}
In this paper we deal with maps of the form $f_\lambda:z\mapsto
\lambda (e^z-1)$, for $\lambda \ge 1$. As long as
$\lambda$ is strictly greater than $1$, $0$ is a repelling fixed point
and there exists an attracting fixed point $q_\lambda<0$.
Those two points collapse to $0$ for $\lambda=1$, and $0$ becomes
parabolic. We are interested in $J_\lambda$, the set of points
that do not escape to $\infty$ under iterations of $f_\lambda$. The
Hausdorff dimension of this set, that we denote $d(\lambda)$,
is an element of $(1,2)$, and is called the Hyperbolic Dimension of
the map $f_\lambda$. While for any $\lambda$ the Julia set of $f_\lambda$ has Hausdorff dimension
constant equal to $2$, cf. \cite{mcmu3}, the Hyperbolic Dimension varies with $\lambda$. Moreover, any
invariant probability measure gives full mass to  $J_\lambda$, and $d(\lambda)$ is, in the hyperbolic case,
equal to the first zero of the pressure of the map $t\mapsto -t\log |f_\lambda'|$\footnote{This result is known as  {\it Bowen's formula}.}, cf. \cite{mu2}.

Variations of $\lambda\mapsto d(\lambda)$ with respect to
$\lambda$, is an interesting feature that reflects changes in geometry
after perturbation of a dynamical system. The philosophy is
that $d$ behaves smoothly, and even real analytically, if we
perturb a conformal hyperbolic dynamical system, in a real analytic way.

This philosophy was proposed in 1981 Rio de Janeiro's conference by Sullivan \cite{su}.
The same year Ruelle  \cite{ru} proved that it was true for a class of Hyperbolic Conformal Repellers. His strategy,
used since then in other contexts, see \cite{uz2} for the exponential family and \cite{mu1} for meromoprhic functions,
was the following : prove a Bowen's formula that identifies the dimension as the zero of a pressure function, prove that
this pressure is the logarithm of a simple and isolated eigenvalue of a Perron-Frobenius(-Ruelle) operator, then use some
results about perturbation theory of operators.

When approaching the boundary of an Hyperbolic components one can not expect any smoothness. Nevertheless there still exists some
paths along thus we still have continuity of the Hausdorff dimension. This was first proved by Bodart and Zinsmeister in \cite{bozi}
for the quadratic family, $z\mapsto z^2+c$, for $c\in\R$ approaching $\frac{1}{4}$ from the left. Then it has been proved for other
parameters $c$, \cite{ri}, or other rational maps, \cite{mcmu2}, \cite{bulei}, or in other situations see \cite{mcmu1} for Kleinian Groups, \cite{uz1}
for the exponential family. The strategy for such results is to control conformal measures, or Patterson-Sullivan measures, in order to
prove that they converge towards the "good" conformal/Patterson-Sullivan-measure. This usually boils down in proving that any limiting
measure is non-atomic. Note that this strategy may also be used to proved discontinuity of the Hausdorff dimension, or more precisely to prove
convergence towards something bigger than the Hausdorff dimension of the "limit set", \cite{dsz}, \cite{uzins1} and \cite{uzins2}.

The problem of the derivative of the Hausdorff dimension is, to our knowledge, investigated in two other papers than the present one.
In \cite{hz} for the quadratic family it is proved that $d'(c)$, the derivative of $d(c):=\mbox{Hdim}(J_c)$, diverges towards $+\infty$
as $c$ converges towards $\frac{1}{4}$ from the left. In \cite{ja}, still for the quadratic family, but this time for $c$ converging from
the right towards $-\frac{3}{4}$, and under the realistic hypothesis that $d(-\frac{3}{4})<\frac{4}{3}$, it is proved that $d'(c)$ converges
towards $-\infty$.
In order to control the derivative the starting point in all those papers is first to get an exact formula for the derivative. This is done using thermodynamic formalism
by differentiating the Bowen's formula. Then some uniform estimates of distorsion in a neighborhood of the fixed point are used in order to
control measures of fondamental annuli. Conclusions then comes from a precise analysis of a certain integral. This is that last point
that explains why such a study has not been yet done in a more general setting. In the present paper, as well as in \cite{hz} and \cite{ja},
some very particular properties of the case studied are used to conclude.

\subsection{Main result}
When one notes that if $\tau_\lambda$ denotes the translation by
$-\lambda$, then we have $f_\lambda\circ\tau_\lambda=\tau_\lambda\circ
g_\lambda$, with $g_\lambda(z)=\alpha(\lambda)e^z$ and
$\alpha(\lambda)=\lambda e^{-\lambda} $, this philosophy (real
analyticity of $d$) is in
\cite{uz2}  proved to be the case. More precisely, it is proved there
that $d:\lambda\mapsto d(\lambda)$ is real-analytic on $(1,+\infty)$,
and in \cite{uz1}, that it is continuous on $[1,+\infty)$.
In this paper we study the asymptotic behavior of the function
$\lambda\mapsto d'(\lambda)$, and we prove the following.

\begin{Th}\label{main}
There exist $\lambda_0>1$ and $K>1$ such that $\forall \lambda\in (1,\lambda_0)$
$$
\left\{
\begin{array}{lclclc}
\frac{-1}{K}(\lambda-1)^{2d(1)-2}
& \le  & d'(\lambda)
& \le & -K(\lambda-1)^{2d(1)-2} & \mbox{if $d(1)<\frac{3}{2}$,}\\\hspace{0.1cm}
& & |d'(\lambda)| &
\le & K (\lambda-1)\log \frac{1}{\lambda-1} & \mbox{if $d(1)=\frac{3}{2}$,}\\\hspace{0.1cm}
& & |d'(\lambda)|
& \le & K(\lambda-1) & \mbox{if $d(1)>\frac{3}{2}$.}
\end{array}
\right.
$$
In particular the function $\lambda\mapsto d(\lambda)$ is C$^1$ on $[1,+\infty)$, with $d'(1)=0$.
\end{Th}

{\bf Remark : }
As already mentioned, conjugating $f_\lambda$ by the translation $\tau_\lambda$, we get the family $g_\lambda:=\tau_\lambda\circ f_\lambda\circ \tau_\lambda^{-1}$,
with $g_\lambda(z)=\lambda e^{-\lambda}e^z$. Changing variable to $\varepsilon:=\lambda e^{\lambda-1}-1$, we get
the family  $g_\varepsilon:z\mapsto (1+\varepsilon)e^{-1}e^z$ with $\varepsilon\sim (\lambda-1)$.
 Let
$D(\varepsilon)$ be the hyperbolic dimension of $g_\varepsilon$,
then
$$\left\{
\begin{array}{rcll}
D'(\varepsilon) & \sim & \varepsilon^{2D(0)-3} & \mbox{if $D(0)<\frac{3}{2}$,}\\\hspace{0.1cm}
|D'(\varepsilon)| & \lesssim & \log \frac{1}{\varepsilon} & \mbox{if $D(0)=\frac{3}{2}$,}\\\hspace{0.1cm}
|D'(\varepsilon)| & \lesssim & K & \mbox{if $D(0)>\frac{3}{2}$.}
\end{array}
\right.
$$
Note in particular that, in case $D(0)<\frac{3}{2}$, we get  exactly the same asymptotic as the one
in \cite{hz} for the family $c\mapsto z^2+c$, with $c<\frac{1}{4}$. For this last family we were able to
prove that $d(\frac{1}{4})<\frac{3}{2}$, see \cite{hz2}. Inequality that we do not know for the exponential family.

Note also that if $d(1)<\frac{3}{2}$ then we have a control on the sign of the derivative in a right neighborhood of $1$. It asserts
that $d(1^+)$ is a local maximum of the Hyperbolic Dimension.
\vskip0.2cm
\noindent
{\bf Remark : }
There is to our knowledge no algorithm to compute accurately Hausdorff dimension of parabolic Julia sets. In \cite{hz2}
an estimate of the Hausdorff dimension of the cauliflower (the Julia set of $z\mapsto z^2+\frac{1}{4}$)
is given using by calculating, with a computer, the first terms of a sum, then by estimating its tail. This method
uses strongly particular properties of the map. More generally,
one could build an infinite iterated function system whose limit set would
have Hausdorff dimension equal to the hyperbolic dimension of the Julia
set. Then, using  results from \cite{hu}, one could
approximate this Hausdorff dimension by finite subsystems keeping track of
the error. Finally, there are algorithms to calculate Hausdorff dimension of finite IFSs
with any desired accuracy \cite{mcmu4}, \cite{jp}. However, to realize such program would be a
tedious extremely time consuming task.

\vskip0.2cm
The proof of the main result will follow exactly the same lines as the
one of \cite{hz2}, but will make an extensive use of
the Thermodynamic Formalism for Meromorphic Functions, as developed
by, Urba\'nski, Urba\'nski and Kotus, Urba\'nski and Zdunik,
and Urba\'nski and Mayer. The reader will find in \cite{mu1} all
proofs of results we need in this paper, as well as a complete
bibliography on the subject.
\subsection{Organization of the paper}
In the first part we use Chapter~8 of \cite{mu1} to  get a formula for
$d'(\lambda)$, for any $\lambda\in(1,+\infty)$. This mainly
consists of conjugating the dynamics and differentiating the pressure.

In the second part we collect some estimates of the distortion around
the fixed point $0$. They are crucial since the formula
obtained in the first part of this paper involves two integrals with
respect to an invariant measure that has unbounded Radon-Nikodym
derivative with respect to the Hausdorff measure, in any neighborhood of $0$.

In the third part we use those estimates to control the integrals and to prove the main result.

In the first appendix we prove the estimates used in the second part
of this paper in a more general setting than needed in this
paper. Namely,
we allow the repelling fixed point to converge towards a parabolic
fixed point with several petals. The second appendix is devoted to the
study of partial sums of some sequences that will be needed several times in the paper.
\vskip0,5cm
{\bf Thanks : } 
The authors thank the european Marie Curie network CODY which help them to meet several times.
They also thank the referee for his suggestions and his careful reading of the paper.

\section{A formula for the derivative of the function $\lambda\mapsto d(\lambda)$}

Before giving and proving the formula of the derivative, this is done below in Proposition \ref{formuladerivative},
we introduce some notations and recall some results concerning the thermodynamical formalism
for that family of exponential maps.

\subsection{Thermodynamic formalism}
Let $P$ be be the cylinder $\{z\in\C\,|\, -\pi<\im z<\pi\}$.
As it is done in \cite{uz2} we associate to $f_\lambda$ the map $F_\lambda:P\to P$ defined by
$$F_\lambda\circ \pi=\pi\circ f_\lambda,$$
with $\pi$ being the natural projection on the cylinder $P=\C/\sim$,
with $z_1\sim z_2$ if and only if $(z_1-z_2)=2ik\pi$,
for some $k\in \Z$. In particular for any $z\in P$ we have
$f_\lambda(z)=F_\lambda(z)$, and $F_\lambda(z)=F_\lambda(z')$ if and
only if there exists $k\in\Z$ such that
$f_\lambda(z)-f_\lambda(z')=2ik\pi$. This
tells us that for any $z\in P$, we have
$F_\lambda^{-1}(z)=\{z_k\in P\,|\ f_\lambda(z_k)=z+2ik\pi,\,
k\in\Z\}$. We also see that
$J(F_\lambda)=\pi(J(f_\lambda))=J(f_\lambda)\cap P$.

Let us now introduce some notation and collect some results, where we
mainly refer to \cite{mu1}, see also \cite{uz2}, \cite{uz1}, \cite{utf}.
\vskip0.2cm\noindent
- For any $\lambda\ge 1$ we define $\l_{\lambda,t}$, the Perron-Frobenius operator associated with the potential
$-t\log |F_\lambda'|$. It acts on $\H_\alpha^\lambda$, the set of bounded
$\alpha$-H\"{o}lder functions defined on $J(F_\lambda)$, in the
following way, let
$g\in \H_\alpha^\lambda$, and $z\in J(F_\lambda)$
$$\begin{array}{rclr}
\l_{\lambda,t}(g)(z) & = &  \ds\sum_{F_{\lambda}(y)=z} \frac{1}{|(F_\lambda)'(y)|^{t}}g(y)& \\ \hspace{0.2cm}
 & = & \ds\sum_{k\in \Z}\frac{1}{|z+\lambda+2ik\pi|^t}g(z_k), & \mbox{with $z_k\in P$ such that $f_\lambda(z_k)=z+2ik\pi\cdot$}
\end{array}$$
\vskip0.2cm\noindent
- The only $d(\lambda)$-conformal measure supported on $J_\lambda$ is denoted $m_\lambda$\footnote{We refer to section 3 of this paper for a definition and more details
about conformal measures.}.
\vskip0.2cm\noindent
- The only equilibrium measure for the potential $-d(\lambda)\log|F_\lambda' |$ and the
dynamical system $(J_\lambda,F_\lambda)$ is denoted $\mu_\lambda$.
\vskip0.2cm\noindent
-The pressure of the potential $-t\log |F_\lambda'|$ is denoted $P(\lambda,t)$, and is defined by
$$
P(\lambda,t)=\sup \{h_\mu-t\chi_\mu\},
$$
where the supremum is taken over all invariant probability measures
$\mu$ supported on $J(F_\lambda)$, such that $\chi_\mu<+\infty$, where
$h_\mu$ denotes the metric entropy of the measure $\mu$, and
$\chi_\mu=\int\log |F_\lambda'|d\mu$ is its Lyapunov exponent.
\vskip0.2cm
We will derive our formula for $d'(\lambda)$ starting from Bowen's
formula that asserts that for any $\lambda>1$, $d(\lambda)$
is the only real number so that $P(\lambda,d(\lambda))=0$ (see \cite{uz2}). We want to
differentiate this formula with respect to $\lambda$, and in order
to do so we need  to appropriately conjugate the dynamics of $F_\lambda$.
\vskip0.2cm\noindent - Let $\lambda_0>1$ be fixed. For any $\lambda
>1$,  we denote  $h_\lambda$ the conjugating map from
$J_{\lambda_0}$ to $J_\lambda$ such that $F_\lambda\circ h_\lambda=h_\lambda\circ F_{\lambda_0}$.
\vskip0.2cm\noindent - We then set : $\varphi_{\lambda,t}:=-t\log |F'_\lambda\circ h_\lambda|$. It is a potential which is defined
on $J_{\lambda_0}$. We then use Corollary 8.10 in \cite{mu1} that tells us that $(\lambda,t)\mapsto P_0(\varphi_{\lambda,t})$ is real
analytic for $\lambda$ close enough to $\lambda_0$\footnote{We denote here $P_0$ the pressure with respect to the dynamical system
$(J_{\lambda_0}, F_{\lambda_0})$.}. Bowen's formula then implies that
$\frac{\partial}{\partial\lambda}P_0(\varphi_{\lambda,d(\lambda)})=0$.
It is this calculation that leads to the desired formula.

\subsection{The formula and its proof}

In this section we prove the following formula
\begin{Prop}\label{formuladerivative}
For any $\lambda\in (1,+\infty)$ we have
\begin{equation}\label{formula}
d'(\lambda)=-\frac{d(\lambda)}{\chi_{\mu_\lambda}}\left(1-\frac{1}{\lambda}\right)\int_{J_\lambda}
\sum_{k=1}^{+\infty} \re\left(\frac{1}{(F_\lambda^k)'}\right)d\mu_\lambda\cdot,
\end{equation}
where $\mu_\lambda$ is the only equilibrium measure for the potential $-d(\lambda)\log|F_\lambda'|$.
\end{Prop}

Let $\lambda_0>1$ be fixed and let $h_\lambda$ denote the conjugating
map : $F_\lambda\circ h_\lambda=h_\lambda\circ F_{\lambda_0}$.
Since $\mu_\lambda$ is the equilibrium measure for the potential $-d(\lambda)\log |F_\lambda'|$, we deduce that
the potential $\varphi_{\lambda,d(\lambda)}$ has a unique equilibrium
measure which is $\tilde{\mu_\lambda}:={h_{\lambda}}_*(\mu_\lambda)$.
We shall now use Theorem 6.14 in \cite{mu1} which asserts that given a
tame function $\varphi$ and a weakly tame function $\psi$ we
have
\begin{eqnarray}\label{derivPressure}
\frac{\partial}{\partial t}P_0(\varphi+t\psi)_{|\, t=0}=\int\psi d\mu_\varphi,
\end{eqnarray}
with $\mu_\varphi$ the equilibrium measure for the potential $\varphi$.
We refer to chapter 4 of \cite{mu1} for definition of tame and loosely tame functions.
By Lemma 8.9 in \cite{mu1}, we know that  for $R>0$ small enough,
there exists $\beta>0$ such that $\forall \lambda\in
(\lambda_0-R,\lambda_0+R)$
$\varphi_{\lambda,t}$ is $\beta$-tame. We then deduce from (\ref{derivPressure}) that
\begin{eqnarray}\label{derivPressure2}
0=\frac{\partial}{\partial \lambda}P_0(\varphi_{\lambda,d(\lambda)})
=\int_{J_{\lambda_0}}\frac{\partial}{\partial \lambda}
\left(\varphi_{\lambda,d(\lambda)}\right) d\tilde{\mu_\lambda}.
\end{eqnarray}
We thus have to compute $\frac{\partial}{\partial \lambda}\varphi_{\lambda,d(\lambda)}$. Note that
$$\varphi_{\lambda,d(\lambda)}=-d(\lambda)\log|F_\lambda'\circ h_\lambda|=-d(\lambda)(\log \lambda+\re h_\lambda).$$
Differentiating with respect to $\lambda$  we get
\begin{eqnarray}\label{derivPressure3}
\frac{\partial}{\partial \lambda}\varphi_{\lambda,d(\lambda)}=-d'(\lambda)\log|F_\lambda'\circ h_\lambda|-d(\lambda)
\left(\frac{1}{\lambda}+\re \frac{\partial}{\partial \lambda}h_\lambda\right)\cdot
\end{eqnarray}

\begin{Lemme}\label{derive1}
For any $\lambda\in (1,+\infty)$ and any $z\in J_{\lambda_0}$ we have
\begin{equation}\label{derive11}
\frac{\partial}{\partial \lambda} h_{\lambda}(z)=\left(1-\frac{1}{\lambda}\right)\sum_{k=1}^{+\infty}
\frac{1}{(F_\lambda^k)'(h_\lambda(z))}-\frac{1}{\lambda}\cdot
\end{equation}
\end{Lemme}
In order to prove this formula we use two results from \cite{utf},
Lemma 13.2 and Proposition 13.4, that we give in the following Lemma
\begin{Lemme}\label{facts}
For any $\lambda_0\in (1,+\infty)$ one can find $R>0$, $K>0$, and $\alpha>0$
such that
\begin{equation}\label{fact1}
\forall \lambda\in B(\lambda_0,R)\quad\forall n\in \N\quad\forall z\in J_{\lambda}\qquad\qquad
|(F_{\lambda}^n)'(z)|\ge K(1+\alpha)^n\cdot
\end{equation}
\begin{equation}\label{fact2}
\forall \lambda\in B(\lambda_0,R)\quad\forall z\in J_{\lambda}\qquad\qquad
\big|\frac{\partial}{\partial \lambda}{h_{\lambda}}(z)\big|  <K\cdot
\end{equation}
\end{Lemme}
We can now prove Lemma \ref{derive1}.
\proof In order to simplify notation, we write ${\dot{h}}_\lambda$
instead of $\frac{\partial}{\partial \lambda}h_\lambda$, and we drop $z$.
We start with the conjugating formula : $h_\lambda\circ F_{\lambda_0}=F_\lambda\circ h_\lambda=\lambda
(e^{h_\lambda}-1)$, that we differentiate with respect to $\lambda$. We thus get,
$$\dot{h_{\lambda}}\circ F_{\lambda_0}=\dot{F_\lambda}\circ h_\lambda+\dot{h_\lambda}F_\lambda'\circ h_\lambda .$$ So that we have
$$\dot{h_{\lambda}}=\frac{\dot{h_{\lambda}}\circ F_{\lambda_0}}{F_\lambda'\circ h_\lambda}-
\frac{\dot{F_\lambda}\circ h_\lambda}{F_\lambda'\circ h_\lambda}\cdot$$
Iterating this formula we end up for $n\in\N$ with
$$\dot{h_{\lambda}}(z)=\frac{\dot{h_{\lambda}}\circ F_{\lambda_0}^n}{(F_\lambda^n)'\circ h_\lambda}-
\sum_{k=1}^n\frac{\dot{F_\lambda}\circ F_\lambda^{k-1}\circ h_\lambda}{(F_\lambda^k)'\circ h_\lambda}\cdot
$$
Using Lemma \ref{facts} we deduce that
$$\frac{\dot{h_{\lambda}}(F_{\lambda_0}^n(z))}{(F_\lambda^n)'(h_\lambda(z))}\quad\mbox{is converging towards 0}\cdot$$
On the other hand, since $\dot{F_\lambda}(z)=e^z-1=\frac{1}{\lambda}F_\lambda'(z)-1$, for any $k$ we have
$$\frac{\dot{F_\lambda}\circ
  F_\lambda^{k-1}}{(F_\lambda^k)'}=\frac{1}{\lambda}\frac{1}{(F_\lambda^{k-1})'}-\frac{1}{(F_\lambda^{k})'}\cdot$$
This leads to
$$\sum_{k=1}^n\frac{\dot{F_\lambda}\circ F_\lambda^{k-1}}{(F_\lambda^k)'}=\frac{1}{\lambda}-\frac{1}{(F_\lambda^n)'}+
\left(\frac{1}{\lambda}-1\right)\sum_{k=1}^{n-1}\frac{1}{(F_\lambda^k)'}\cdot$$
Using (\ref{fact1}) in Lemma \ref{facts} we get that the series on the left above is converging towards
$$\sum_{k=1}^{+\infty}\frac{1}{(F_\lambda^k)'},$$
which finishes the proof.
\cqfd

Using (\ref{derivPressure3}) and Lemma \ref{derive1} in (\ref{derivPressure2}) we get
\begin{eqnarray}\label{derivPressure4}
\end{eqnarray}
$\displaystyle -d'(\lambda)\int_{J_{\lambda_0}} \log |F_\lambda'| d\tilde{\mu}_\lambda-d(\lambda)\left(1-\frac{1}{\lambda}\right)
\int_{J_{\lambda_0}} \re\sum_{k\ge 1}\frac{1}{(F_\lambda^k)'\circ h_\lambda}=0\cdot$
\\For any function $g$ continuous on $J_\lambda$ we have $\tilde{\mu}_\lambda(g\circ h_\lambda)=\mu_\lambda(g)$.
We deduce from (\ref{derivPressure4}) that Proposition \ref{formuladerivative} is true.

\section{Local dynamic and uniform estimates}
In this section we introduce some notations and collect estimates
proved in the appendix in a more general setting\footnote{We deal
in the appendix with a family of germ of holomorphic in a neighborhood
of a repelling fixed point which degenerates into a parabolic fixed
point
with $p$ petals.}. We then use these estimates in order to control
uniformly conformal measures $(m_\lambda)$ and equilibrium
measures $(\mu_\lambda)$.
\subsection{Notation}

We know that
$J_\lambda\cap P\subset\{z\in\C\,|\, -\frac{\pi}{2}<\im z<\frac{\pi}{2}\}$.
\vskip0.2cm\noindent
Given $0<\theta<\frac{\pi}{2}$ we denote $S_\theta$
the sector $\{re^{i\alpha}\,|\, r>0,\,-\theta<\alpha<\theta\,\}$.
\vskip0.2cm\noindent For $r_0<<1$ we fix $0<\theta<\frac{\pi}{2}$ to be such that $J_1\cap B(0,r)\subset S_\theta$. Then we choose $\varepsilon_0>0$ small enough
so that for any $0\le\lambda=1+\varepsilon\le \lambda_0=1+\varepsilon_0$ we have $f_\lambda^{-1}(S_\theta)\subset S_\theta$ and  $J_\lambda\cap B(0,r)\subset S_\theta$.
We then set $\gamma_0=\{r_0e^{it}\,|\ t\in]-\theta,\theta[\}$, $\gamma_1(\lambda)=f_\lambda^{-1}(\gamma_0)$. Joining $r_0e^{i\theta}$ with
$f_\lambda^{-1}(r_0e^{i\theta})$ by a line, and doing the same with $r_0e^{-i\theta}$ and its image by $f_\lambda^{-1}$, we get a cell $C_0(\lambda)$. It is
a simply connected domain. A compactness argument tells us that if $1\le\lambda\le\lambda_0$, then there exists a simply connected domain $V\subset S_\theta$ such that
the closure of $\cup_{\lambda}C_0(\lambda)$ is a subset of $V$. In particular, K{\oe}be distorsion Theorem gives us a constant $K>1$, only depending on $r_0$ and $\lambda_0$, such that for any univalent function $h$ on $V$ and any point $x$ and $y$ in $\overline{\cup_\lambda C_0(\lambda)}$ we have $\frac{1}{K}\le \frac{|h'(x)|}{|h'(y)|}\le K$.
We will use later on this fact with inverse branches of $f_\lambda^n$. They are well defined on $V$ since the post-singular set of the $f_\lambda$'s, i.e. the orbit $-\lambda$
under $f_\lambda$, is a subset of $(-\infty,0)$.

We then define for each integer $n$ the set
$C_n(\lambda):=f_\lambda^{-n}(C_0(\lambda))$, with $f_\lambda^{-n}$
being  the $n^{\mbox{\footnotesize th}}$ iterates of the inverse branch of $f_\lambda$ defined on
$B(0,r_0)$ that fixes $0$. In the following we are working with
respect to measures concentrated on $J_\lambda$ of dimension strictly
greater to $1$. One checks easily
in that context that with respect to such measure
$(C_n(\lambda))_{n\in\N\cup\{0\}}$ is a partition of
$B(0,r_0)$. Moreover the set $C_0(\lambda)$ is mapped univalently
by $f_\lambda^{-n}$ to $C_n(\lambda)$.
\vskip0.2cm\noindent
Let $N_\varepsilon$ be an integer\footnote{In our study we have $N_\varepsilon\sim \frac{1}{\varepsilon}=\frac{1}{\lambda-1}\cdot$} and defined the sequence
$(a_{n,\varepsilon})_{n\in\N}$ as $a_{n,\varepsilon}=\frac{1}{n}$, if $n\le N_\varepsilon$, and $a_{n,\varepsilon}=\varepsilon (1+\varepsilon)^{-n}$, if
$n\ge N_\varepsilon$. Note that $a_{n,\varepsilon}\to 0$. In order to simplify notations, we let $a_n:=a_{n,\varepsilon}$. We consider now the one parameter familly
of sequences, $(a_n(\alpha))_{n\in\N}$, defined for $n\in\N$ by $a_n(\alpha):=a_n^\alpha$. We are also interested in partial sums of $\sum a_n(\alpha)$. For
$k\le n$ we let $S_{k,n}(\alpha):=\sum_{l=k}^n a_l(\alpha)$. The sequence $(a_n(\alpha))$ will describe, for different values of $\alpha$, the distorsion
around $0$, the conformal measure of partition sets of a neighborhood of $0$, and the partial sums $S_{k,n}(\alpha)$ will play a role in controlling the
invariant measure of the same partition sets, as well as evaluating the integral which is crucial in order to get our main result. Those estimates
are easy and we use them in this section but we postponed their proofs to the appendix.

\subsection{Uniform estimates of the distorsion}
In this section we give uniform estimates depending on $\lambda$ for
the local dynamics next to the repelling-parabolic fixed point $0$. We recall
that the family we are studying is given for $\lambda:=1+\varepsilon\ge 1$ by
$f_\lambda(z)=\lambda (e^z-1)$. In particular, in a neighborhood of $0$,
the local dynamic is given by the following Taylor expansion
$$F_\lambda(z)=f_\lambda(z)=\lambda z+z^2+z^3g_\lambda(z)\cdot$$
With $g_\lambda(z)$ uniformly bounded, independently of $\lambda$, as
soon as a neighborhood of zero has been fixed. Note in particular that for $\varepsilon=0$
, the point $0$ is a parabolic fixed point with one petal.
\vskip0.2cm\noindent
We apply the general results of the first appendix of this paper to this special family $f_\lambda$.
In the remaining of the paper we set $\lambda=1+\varepsilon$ and we
denote the relevant quantities by indexing them equally well either by $\varepsilon$ or
$\lambda$. Moreover, in the remainder of this section
$F_\lambda^{-n}$ will be the inverse branch of $F_\lambda^{-n}$ that
fixes $0$. From Proposition \ref{unifestimatesC} we deduce that
\begin{Prop}\label{unifestimatesCp}
Let $0<r_0$, $1<\lambda_0$ being fixed. Then there exists $K>1$ such
that $\forall  \lambda\in (1,\lambda_0)$,  $\forall z\in
C_0(\lambda)$, and $\forall n\in\N$
$$\frac{1}{K}a_n(2)\le |(F_\lambda^{-n})'(z)|\le K a_n(2)\cdot$$
\end{Prop}
The following result is technical but will be crucial in order to control the sign of the derivative $d'(\lambda)$.
\begin{Lemma}\label{Argumentznp}
Let $0<r_0$, $1<\lambda_0$ being fixed. There exists an integer $N$
such that $\forall n>N$, $\forall  k\in \N\cap [1,n-N]$,
$\forall \lambda\in (1,\lambda_0)$ and  $\forall z\in C_n(\lambda)$
$$ \frac{\sqrt{3}}{2}|(F_\lambda^k)'(z)|\le \re (F_\lambda^k)'(z)\cdot$$
\end{Lemma}
\proof Let $z\in C_n$ and $\theta_k(z)=\arg (F_\lambda^k)'(z)$. The
Lemma boils down to proving that $|\theta_k(z)|\le \frac{\pi}{6}$.
\\One computes that $(F_\lambda^k)'(z)=\prod_{j=0}^{k-1}F_\lambda'(F_\lambda^j(z))=\lambda ^k\exp (\sum_{j=0}^{k-1} F_\lambda^j(z))$. So that
we have $\theta_k(z)=\sum_{j=0}^{n-1}\im (F_\lambda^j(z))$. Since $F_\lambda^j(z)$ belongs to $C_{n-j}$ we may use Corollary \ref{Argumentzn}
which asserts that $|\im (Z)|\lesssim \frac{1}{(n-j)^2}$ for any $Z\in C_{n-j}$. We thus have

$\ds |\theta_k(z)| \lesssim  \sum_{j=1}^{n-k}\frac{1}{(n-j)^2}\le \sum_{j=N}^{+\infty}\frac{1}{j^2}.$

This is less than $\frac{\pi}{6}$ if $N$ is big enough and we are done.
\cqfd

\

We end this section with two more estimates of the distorsion.
The first one needs the following observation on the localization of $J(f_\lambda)$.

\begin{Lemma}\label{l1062805}
For every $R>0$ there exists $\De>0$ such that for all $\lam>1$,
$$J(f_\lam)\setminus \bigcup_{n=-\infty}^{+\infty} B(2\pi ni,R)
\subset \{z\in\C:\re z\ge \De\}.$$
\end{Lemma}

{\sl Proof.} First notice that
$$f_\lam(\{z\in\C:\re z<0\})=B(-\lam,\lam)\subset \{z\in\C:\re z<1-\lam\}\subset \{z\in\C:\re z<0\}.$$
Thus
\begin{equation}\label{1062805}
\{z\in\C:\re z<0\}\subset {\mathcal F}(f_\lam):=\mbox{Fatou set of $f_\lambda$}.
\end{equation}

Now write $z=x+iy$. Then
$$\re(f_\lam(z))=\re(\lam(e^x\cos y+ie^x\sin y-1))=\lam(e^x\cos y-1).$$

Note that there exists $\De_1>0$ so small that if $0<x<\De_1$
and $x+iy\notin \bigcup_{n=-\infty}^{+\infty} B(2\pi ni,R)$,
then $\dist(y, \{2\pi ni:n\in\Z\})>R/2$, and consequently,
$\cos y<\cos(R/2)$. Hence, $\re(f_\lam(z))<\lam(e^{\De_1}\cos(R/2)-1)$. Take now $0<\De\le\De_1$
so small that $e^{\De}\cos(R/2)<1$. So $\re(f_\lam(z))<0$ and, by
(\ref{1062805}), $f_\lam(z)\in{\mathcal F}(f_\lam)$. Therefore we have
proved that
$$\{z\in\C\setminus \bigcup_{n=-\infty}^{+\infty} B(2\pi ni,R):\re z<\De\}\
\subset {\mathcal F}(f_\lam).$$ We are done.
\cqfd

\noindent Now notice that if $\re z\ge\De$, then
$$|f_\lam'(z)|=\lam e^{\re z}\ge \lam e^{\De}>1.$$
Combining this and Lemma~\ref{l1062805}, we obtain the following.
\
\begin{Lemma}\label{l2062805}
For every $R>0$ there exists $\gamma>1$ such that for every $z\in
J(F_\lam) \setminus B(0,R)$,
$$
|F_\lam'(z)|\ge\gamma.
$$
\end{Lemma}
Using Proposition \ref{unifestimatesC}, Lemma \ref{l2062805} and the same reasoning as for the proof of Lemma 3.6
in \cite{hz} we prove the following result
\begin{Lemma}\label{landingoutsideB}
There exist $0<r_0$, $1<\lambda_0$ and $1<K$ such that $\forall \lambda\in (1,\lambda_0)$ and $\forall z\in J_\lambda$,
$$F_\lambda^n(z)\notin B(0,r_0)\Rightarrow Kn^2\le |(F_\lambda^n)'(z)|.$$
\end{Lemma}

\subsection{Conformal measures}
Let us recall that a probability
measure $m_\lambda$ is called {\it conformal} if its strong Jacobian
is equal to $|F_\lambda'|^{d(\lambda)}$. This means that for any
measurable set $A$ on which $f_\lambda$ is 1-1 we have
\begin{equation}\label{defconfmeas}
m_\lambda(F_\lambda(A))=\int_A |F_\lambda'|^{d(\lambda)}
dm_\lambda\cdot
\end{equation}
Those measures are usually a powerful tool to study Hausdorff
dimension of Julia sets. In fact their definition is dynamical but
they very often carry a geometrically significant information about the Julia
set. In many of cases they coincide (up to a multiplicative constant)
with Hausdorff or packing measures on the Julia
set.

Using Proposition \ref{unifestimatesC} and the notation introduced
below we get the following.

\begin{Prop}\label{unifestimatesconf}
Let $0<r_0$, $1<\lambda_0$ being fixed. Then there exists $K>1$ such
that $\forall  \lambda\in (1,\lambda_0)$,  and $\forall n\in\N$
$$
 \frac{1}{K} a_n(2d(\lambda)) \le m_\lambda(C_n(\lambda)) \le K a_n(2d(\lambda))
$$
\end{Prop}
\proof This is not difficult when one observes that for each $\lambda$
the function $F_\lambda^n$ is univalent on $C_n(\lambda)$. In
particular using the definition of a conformal measure we deduce that
:
\begin{equation}\label{B1}
m_\lambda (C_0(\lambda))=m_\lambda(F_\lambda^n(C_n(\lambda)))=\int_{C_n(\lambda)}
|(F_\lambda^n)'|^{d(\lambda)}dm_\lambda .
\end{equation}
We then use estimates of Proposition \ref{unifestimatesCp}, since
$|(F_\lambda^n)'|$ on $C_n(\lambda)$ is comparable with
$|(F_\lambda^{-n})'|^{-1}$  on  $C_0(\lambda)$. We deduce that there exists a
constant $K>0$ such that for any $z\in C_0(\lambda)$
\begin{equation}\label{B2}
\frac{1}{K}|(F_\lambda^n)'(z)|^{-d(\lambda)}\le
m_\lambda(C_n(\lambda))\le K |(F_\lambda^n)'(z)|^{-d(\lambda)}\cdot
\end{equation}
We can now conclude the proof by using again Proposition \ref{unifestimatesC}.
\cqfd

\

{\bf Remark :} Let $m_\infty$ be any accumulation point of the
family of probability measures $(m_\lambda)_{\lambda> 1}$. Let
$(\lambda_n)$ be a sequence of real numbers converging from above
towards $1$ such that the sequence $(m_{\lambda_n})$ converges weakly to
$m_\infty$, and $(d(\lambda_n))$ converges to some $d\ge 0$. For
any $r>0$ small enough one may find $N(r)$ such that
$$
\begin{array}{r}
\forall n\ge N(r)\qquad B(0,r)\cap J_{\lambda_n}\subset \{0\}
\cup\bigcup_{k\le N(r)}C_k(\lambda_n)\\
\end{array}
$$
And in particular we conclude if $r>0$ is such that
$m_\infty(\{|z|=r\})=0$, that we have
$$m_\infty(B(0,r))=\lim_{n\to\infty}m_{\lambda_n}(B(0,r))\le\lim_{n\to\infty}
\sum_{k>N(r)} m_{\lambda_n}(C_k(\lambda_n))\le \frac{K}{N(r)^{2d(\lambda_n)-1}}\cdot$$
So that we conclude that $m_\infty$ has no atom at $0$. And it
is one of the main point in order to conclude that $d(\lambda)\to
d(1)$ when $\lambda\to 1$, see \cite{uz1}.

We end this section about conformal measures with a technical Lemma. It will be used in the next
section concerning invariant measures.

Before stating and proving this result we recall that $P=\{z\in\C\, |\, -pi<\im z\le pi\}$, and for any
$M>0$, and any $r>0$, we introduce the following notation : $P_M:=\{z\in  P\,|\, \re z\le M\}$, $\b_r:=P\setminus B(0,r)$
and $\b_{r,M}:=P_M\cap \b_r$.

\begin{Lemme}\label{confmeasineq}
There exists $0<\alpha<\beta$ such that $\forall M\ge 2$,
$\forall\lambda\in [1,\lambda_0]$, with $\lambda_0<\frac{\pi}{3}$,
$\forall r\in ]0,\frac{\pi}{3}-\lambda_0]$ and
$\forall A\subset B(0,r)$ measurable, we have
\begin{equation}\label{confmeasineqa}
\alpha m_\lambda(A)\le m_\lambda(F_\lambda^{-1}(A)\cap\b_r)
\le \beta m_\lambda(A).
\end{equation}
and
\begin{equation}\label{confmeasineqb}
m_\lambda(F_\lambda^{-1}(A)\cap\b_{r,M})\le m_\lambda(F_\lambda^{-1}(A)\cap\b_r)
\le 54\beta m_\lambda(F_\lambda^{-1}(A)\cap\b_{r,M}).
\end{equation}
\end{Lemme}
\proof
Let $B_k$ be the connected component of $F_\lambda^{-1}(B(0,r))$ such that $f_\lambda(B_k)=B(2ik\pi,r)$. For any $z\in B_k$ we have :
\begin{equation}\label{confmeasineq1}
|F_\lambda'(z)|=|f_\lambda'(z)|=|f_\lambda(z)+\lambda|=|\lambda+2ik\pi+a e^{i\theta}|,
\end{equation}
with $a<r$ and $0\le \theta<2\pi$. With our assumptions this leads, for $|k|\ge 1$, to
\begin{equation}\label{confmeasineq2}
2|k|\pi-\frac{\pi}{3}\le |F_\lambda'(z)|=|f_\lambda'(z)|\le 2|k|\pi+\frac{\pi}{3}\cdot
\end{equation}
Since $|f_\lambda'(z)|=\lambda\exp(\re z)$, we also get, for any $|k|\ge 1$, that
$$\forall z\in B_k\quad\log 5\le \re z\cdot$$
As a consequence we see that $F_\lambda^{-1}(B(0,r))\cap \b_r=\cup_{|k|\ge 1}B_k$.

The measure $m_\lambda$ being conformal we have
$$m_\lambda(A)=m_\lambda(F_\lambda(A_k))=\int_{A_k}|F_\lambda'|^{d(\lambda)}dm\lambda,$$
with $A_k:=F_\lambda^{-1}(A)\cap B_k$. From (\ref{confmeasineq2}) we deduce that
\begin{equation}\label{confmeasineq3}
\frac{m_\lambda(A)}{(2|k|\pi+\frac{\pi}{3})^{d(\lambda)}}\le m_\lambda(A_k)\le \frac{m_\lambda(A)}{(2|k|\pi-\frac{\pi}{3})^{d(\lambda)}},
\end{equation}
so that
\begin{equation}\label{confmeasineq4}
2 m_\lambda(A)\sum_{k\ge 1}\frac{1}{(2k\pi+\frac{\pi}{3})^{d(\lambda)}}\le m_\lambda(F_\lambda^{-1}(A)\cap \b_r)\le
2 m_\lambda(A)\sum_{k\ge 1}\frac{1}{(2k\pi-\frac{\pi}{3})^{d(\lambda)}}\cdot
\end{equation}
The function $\lambda\mapsto d(\lambda)$ being continuous on $[1,\lambda_0]$ one may consider its minimum $\delta$ which is strictly greater than
1. With $\alpha= 2\sum_{k\ge 1}\frac{1}{(2k\pi+\frac{\pi}{3})^{2}}$ and $\beta=2\sum_{k\ge 1}\frac{1}{(2k\pi-\frac{\pi}{3}){\delta}}$ we have :
$$\alpha m_\lambda(A)\le m_\lambda(F_\lambda^{-1}(A)\cap \b_r)\le \beta m_\lambda (A).$$
This is (\ref{confmeasineqa}).

Note that (\ref{confmeasineq2}) tells us that for any $z\in B_1$ we have $\re z\le\log (2\pi+\frac{\pi}{3})<2\le M$. This implies that
$B_1\subset \b_{r,M}$. In particular we have
$$m_\lambda(A_1)\le m_\lambda(F_\lambda^{-1}(A)\cap \b_{r,M})\cdot$$
We then deduce from (\ref{confmeasineq3}) that
$$\frac{m_\lambda(A)}{(2\pi+\frac{\pi}{3})^2}\le m_\lambda(F_\lambda^{-1}(A)\cap \b_{r,M})\cdot$$
Together with (\ref{confmeasineqa}) we conclude that
$$m_\lambda(F_\lambda^{-1}(A)\cap \b_r)\le (2\pi+\frac{\pi}{3})^2 \beta m_\lambda(F_\lambda^{-1}(A)\cap \b_{r,M})\cdot$$
Since $(2\pi+\frac{\pi}{3})^2\le 54$ we conclude that the left hand side inequality of (\ref{confmeasineqb}) holds.
The right hand side being obvious the proof is finished.
\cqfd

\subsection{Invariant measures}
Let us first recall that $\mu_\lambda=\rho_\lambda m_\lambda$ is the
unique $F_\lambda$-invariant probability measure equivalent with $m_\lambda$.
This measure is also the unique equilibrium state for the potential $-d(\lambda)\log |F_\lambda'|$ i.e.
$$
h_{\mu_\lambda}-d(\lambda)\int\log|F_\lambda'|d\mu_\lambda=\sup\{h_{\mu}-d(\lambda)\int\log|F_\lambda'|d\mu\},
$$
where supremum is taken over all $F_\lambda$-invariant ergodic
probability measures such that $\int\log|F_\lambda'|d\mu<+\infty$. The function $\rho_\lambda$ is obtained
in \cite{mu1} as the limit of the sequence $\l_\lambda^n(1)$. The main results of this section is
\begin{Prop}\label{invariantmeasure}
Let $0<r_0$, $1<\lambda_0$ being fixed. Then there exists $K>1$ such
that $\forall  \lambda\in (1,\lambda_0)$,  and $\forall n\in\N$
$$
\begin{array}{rrclr}
\mbox{i-} & \frac{1}{K} a_n(2d(\lambda)-1) \le &
  \mu_\lambda(C_n(\lambda))& \le K a_n(2d(\lambda)-1)& \mbox{if
  $n\le N_\varepsilon$.} \vspace{0.2cm} \\
\mbox{ii-} & \frac{1}{K}\frac{a_n(2d(\lambda))}{\varepsilon}  \le &
 \mu_\lambda(C_n(\lambda))  & \le K\frac{a_n(2d(\lambda))}{\varepsilon}&
 \mbox{if $n\ge N_\varepsilon$.} \\
\end{array} $$
\end{Prop}
\proof
Let $\b_r:=P\setminus B(0,r)$. We know that $\mu_\lambda(\b_r)>0$ so that the first return time
$N_{\lambda,r}(z):=\inf\{n\ge 1\,|\,F_\lambda^n(z)\in \b_r\}$ is finite $\mu_\lambda$-almost-surely.
Let $\b_{\lambda,n}:=\{N_{\lambda,r}=n\}$. We recall that the sets $(C_n)$ are introduced at the beginning of this section.
Note that for $r$ small enough we have $\b_{\lambda,n}\cap B(0,r)=C_{n-1}(\lambda)$.
Since $\mu_\lambda$ is $F_\lambda$-invariant its restriction to $\b_r$ is invariant for the first return map in $\b_r$,
that we denote $T_\lambda$. Moreover, $\mu_\lambda$ can be built from this $T_\lambda$-invariant measure
and this leads, for any measurable set $A$, to the formula
$$\mu_\lambda(A)=\sum_{n\ge 1}\sum_{k=0}^{n-1}\mu_\lambda(F_\lambda^{-k}(A)\cap\b_{\lambda,n}\cap\b_r).$$
We are interested in the sets $C_l$ for which we get
$$
\mu_\lambda(C_l)=\sum_{n\ge
  1}\sum_{k=0}^{n-1}\mu_\lambda(F_\lambda^{-k}(C_l)\cap\b_{\lambda,n}\cap\b_r).
$$
Note now that the set $F_\lambda^{-k}(C_l)\cap\b_{\lambda,n}\cap\b_r$ is empty
unless $n>l+1$ and $k=n-l-1$. In this case
we have $F_\lambda^{-(n-l)}(C_l)\cap\b_{\lambda,n}\cap\b_r=F_\lambda^{-1}(C_{n-2})\cap\b_r$. We thus conclude that
$$\mu_\lambda(C_l)=\sum_{n\ge l}\mu_\lambda(F_\lambda^{-1}(C_n)\cap\b_r)\cdot$$
In Corollary \ref{ineqconfinv1}, that we admit for the moment, we show that there exists $K_1>0$, independent of
$\lambda$, such that for any $A\subset B(0,r)$ we have
$$\frac{1}{K_1}m_\lambda(A)\le\mu_\lambda(F_\lambda^{-1}(A)\cap\b_r)\le K_1 m_\lambda(A)\cdot$$
So,
$$\frac{1}{K_1}\sum_{n\ge l}m_\lambda(C_n)\le\mu_\lambda (C_l)\le K_1\sum_{n\ge l}m_\lambda(C_n).$$
From Proposition~\ref{unifestimatesconf} we deduce that there exists $K_2>0$ such that
$$\frac{1}{K_2}\sum_{n\ge l}\sum_{n\ge l}a_n(2d(\lambda))\le \mu_\lambda (C_l)\le K_2\sum_{n\ge l}a_n(2d(\lambda)).$$
With the notations used in the appendix this is exactly
$$\frac{1}{K_2} S_{l,+\infty}(2d(\lambda))\le\mu_\lambda (C_l)\le K_2 S_{l,+\infty}(2d(\lambda))\cdot$$
We then use Corollary \ref{partialsum2} to finish the proof.
\cqfd
\begin{Lemma}\label{ineqconfinv}
There exists $K>0$ such that for all $\lambda=1+\varepsilon$, with $\varepsilon>0$ small enough, all $r>0$ small
enough, and all $M>0$ big enough we have,
$$\frac{1}{K}
\le \rho_\lambda
\le K\quad\mbox{on $\b_{r,M}$},\quad\mbox{and}\quad\rho_\lambda
\le K\quad\mbox{on $\b_r$}\cdot
$$
\end{Lemma}
From this Lemma and Lemma \ref{confmeasineq} we easily conclude this.
\begin{Cor}\label{ineqconfinv1}
There exists $K>0$ such that for all $\lambda=1+\varepsilon$, with $\varepsilon>0$ small enough,
all $r>0$ small enough, and for any measurable set $A\subset B(0,r)$ we have
$$\frac{1}{K}m_\lambda(A)\le\mu_\lambda(F_\lambda^{-1}(A)\cap\b_r)\le Km_\lambda(A)\cdot$$
\end{Cor}
\proof
Let $r>0$ and $\varepsilon>0$ be small enough so that the assertions of
Lemmas~\ref{confmeasineq} and \ref{ineqconfinv}
hold. Let $K>0$ coming from Lemma~\ref{ineqconfinv} be larger than
$\max\{\beta,\alpha^{-1}\}$, both $\alpha$ and $\beta$ coming from
Lemma~\ref{confmeasineq}.
By Lemma \ref{confmeasineq} we know that for any $A\subset B(0,r)$ we have
$$\frac{1}{K}m_\lambda(F_\lambda^{-1}(A)\cap\b_r)\le m_\lambda (A)\le K m_\lambda(F_\lambda^{-1}(A)\cap\b_r)\cdot$$
From the right hand side inequality in Lemma \ref{ineqconfinv} we know that
$$\mu_\lambda (F_\lambda^{-1}(A)\cap\b_r)\le K m_\lambda(F_\lambda^{-1}(A)\cap\b_r)\cdot$$
These two inequalities give us
$$\mu_\lambda (F_\lambda^{-1}(A)\cap\b_r)\le K^2 m_\lambda (A)\cdot$$
For the other inequality we first note that Lemma \ref{confmeasineq} also asserts that
$$\frac{1}{K}m_\lambda(F_\lambda^{-1}(A)\cap\b_{r,M})\le m_\lambda(F_\lambda^{-1}(A)\cap\b_r)\le K m_\lambda(F_\lambda^{-1}(A)\cap\b_{r,M})\cdot$$
Since Lemma \ref{ineqconfinv} implies that
$$\frac{1}{K}m_\lambda(F_\lambda^{-1}(A)\cap\b_{r,M})\le \mu_\lambda(F_\lambda^{-1}(A)\cap\b_{r,M})\le K m_\lambda(F_\lambda^{-1}(A)\cap\b_{r,M}),$$
we conclude that
$$m_\lambda (A)\le K m_\lambda(F_\lambda^{-1}(A)\cap\b_r)\le K^2 m_\lambda(F_\lambda^{-1}(A)\cap\b_{r,M})\le K^3
\mu_\lambda(F_\lambda^{-1}(A)\cap\b_{r,M})\cdot$$
We easily deduce that
$$m_\lambda (A)\le K^3 \mu_\lambda(F_\lambda^{-1}(A)\cap\b_r)\cdot$$
This is the left hand side inequality of the Corollary and its proof is finished.
\cqfd

\proof
Before starting the proof of Lemma \ref{ineqconfinv} we sketch the strategy. We first use a result of Urba\'nski and Zdunik, Lemma 3.4
in \cite{uz2}, that asserts that as long as we stay far away from the post-singular set, iterates of $\l_\lambda$ are uniformly
bounded from above by a constant that does not depend on $\lambda$. This gives us that $\rho_\lambda$ is bounded from above in some $\b_r$.
And this allows us to prove that for $r$ and $\varepsilon$ small enough, and for $M$ big enough we have
$$\frac{1}{2}\le \mu_\lambda(\b_{r,M})\le 1.$$

In order to control $\rho_\lambda$ on $\b_{r,M}$ we use K\oe be's distortion Theorem on $\b_{r,M}$ and prove
that the measures $m_\lambda$ have the bounded distortion property on $\b_{r,M}$, with a constant which only depends on $r$ and $M$.
This implies, see \cite{ma} (compare \cite{ha} Propositions 1.2.7 and 1.2.8), that there exists an $F_\lambda$-invariant measure $\nu_\lambda$
which gives mass $1$ to $\b_{r,M}$ and which is equivalent with $m_\lambda$. Its Radon-Nikodym derivative is such that
$\frac{1}{K}\le \frac{d\nu_\lambda}{dm_\lambda}\le K$ on $\b_{r,M}$,
with some $K>0$ independent of $\lambda$. Since $m_\lambda$ is
ergodic and conservative, there is, up to a multiplicative constant,
only one possible invariant measure equivalent to it. This means that
$\mu_\lambda=\alpha_\lambda\nu_\lambda$. Integrating on $\b_{r,M}$ we
conclude that $\alpha_{\lambda}=\mu_\lambda(\b_{r,M})$.
This leads to $\frac{1}{2K}\le \rho_\lambda\le K$.

We now go into further details. Note that the
singular set of $F_\lambda$ is the one point $-\lambda$ which sequence
of iterates converges towards $0$ from the left. In particular
$\b_{r,M}$ is a simply connected domain on which inverse branches of
$F_\lambda$ are well defined. Since $J_\lambda$ is a subset of
$\{-\frac{\pi}{2}\le\im z\le\frac{\pi}{2}\}$ one may find an open
simply connected domain $U_{r,M}$ such that : $\overline{U_{r,M}}\subset
\b_{\frac{r}{2},2M}$
and $J_\lambda\cap \b_{r,M}\subset U_{r,M}$. We have thus an annulus
$\b_{\frac{r}{2},2M}\setminus U_{r,M}$ and an associate K\oe be
constant
$\sqrt{K_{r,M}}$. We conclude that for any $\lambda$ and any $n\in\N$ we have
\begin{equation}\label{ineqconfinv00}
\forall x\in U_{r,M}\quad\forall y\in U_{r,M}\qquad\frac{1}{K_{r,M}}\le \frac{\l_{\lambda}^n(1)(x)}{\l_{\lambda}^n(1)(y)}\le K_{r,M}\cdot
\end{equation}
Since for a measurable set $A$ we have $m_\lambda(F_\lambda^{-n}(A))=\int_A\l_\lambda^n(1)dm_\lambda$, we conclude, if $A\subset U_{r,M}$, that
$$\frac{1}{K_{r,M}}\frac{m_\lambda(A)}{m_\lambda(U_{r,M})}\le \frac{m_\lambda(F_\lambda^{-n}(A))}{m_\lambda(F_\lambda^{-n}(U_{r,M}))}\le
K_{r,M}\frac{m_\lambda(A)}{m_\lambda(U_{r,M})}\cdot$$
This is precisely the bounded distortion property for $m_\lambda$ on
$U_{r,M}$ as it is used in \cite{ha}. Since $(J_\lambda,F_\lambda,m_\lambda)$
is ergodic and conservative there is, up to a multiplicative constant,
only one invariant measure equivalent with $m_\lambda$. Let
$\nu_\lambda$ be the one that gives mass 1 to $\b_{r,M}$. It follows
from Propositions~1.2.7 and 1.2.8 in \cite{ha} that
$$
m_\lambda\mbox{-almost surely on $\b_{r,M}$}
\qquad\frac{1}{K_{r,M}}\le\frac{d\nu_\lambda}{dm_\lambda}\le K_{r,M}\cdot
$$
The measures $\mu_\lambda$ and $\nu_\lambda$ only differ by a
multiplicative constant which can be computed by integrating
the function 1 over $\b_{r,M}$. We deduce that $\mu_\lambda=\mu_\lambda(\b_{r,M})\nu_\lambda$ and we conclude that
\begin{equation}\label{ineqconfinv01}
m_\lambda\mbox{-almost surely on $\b_{r,M}$}\qquad\frac{\mu_\lambda(\b_{r,M})}{K_{r,M}}\le\rho_\lambda\le K_{r,M}
\mu_\lambda(\b_{r,M}).
\end{equation}
Using inequalities (\ref{ineqconfinv00}) one may now adapt the  reasoning of Lemma 3.4 in \cite{uz2} to our situation.
Let $M$ be large enough and $r$ small enough so that : $\frac{\log M-1}{2}\ge r$ and for all $\lambda\in[1,\lambda_0]$
if $\re z>M$ then $\l_\lambda (1)(z)\le 1$. The purpose of the first requirement is the following
\begin{equation}\label{ineqconfinv02}
\forall z\in P\qquad (\re z>M\quad\mbox{and}\quad F_\lambda(y)=z)\Rightarrow |y|>r\quad\mbox{(i.e. $y\in\b_r$)}\cdot
\end{equation}
We prove by induction that $H_n$ is true for all $n$ with
$$H_n\Leftrightarrow||\l_\lambda^n(1)\chi_{B_r}||_\infty\le\frac{K_{r,M}}{m_\lambda(\b_{r,M})}\cdot$$
Notice that $H_0$ is obvious and assume that $H_n$ is true. Since
$$\l_\lambda(1)(z)\le \sum_{k\ge \re z}\frac{2}{k^{d(\lambda)}},$$ and since
$d(\lambda)$ is converging towards $d(\lambda_0)$, one deduces that
$\l_\lambda (1)(z)$ is, uniformly in $\lambda$, converging towards
$0$ as $\re z\to\infty$. We deduce that
$||\l_\lambda(1)\chi_{\b_r}||_\infty$ is achieved for some $z_1\in
\b_r$. An easy induction leads, for
all integers $n\ge 0$, to the existence of some $z_n\in\b_r$ such that
$$
||\l_\lambda^n(1)\chi_{B_r}||_\infty=\l_\lambda^n(1)(z_n)\cdot
$$
Consider $z_{n+1}$ and assume that it lies in $\b_{r,M}$. Then we have
$$
1=\int\l_\lambda^{n+1}(1)dm_\lambda\ge
\int\l_\lambda^{n+1}(1)\chi_{\b_{r,M}}dm_\lambda\ge
\frac{\l_\lambda^{n+1}(1)(z_{n+1})}{K_{r,M}}
m_\lambda(\b_{r,M})\cdot
$$
The last inequality is an application of (\ref{ineqconfinv00}) and we
conclude that $H_{n+1}$ is true. But $z_{n+1}$ might
be with a real part greater than $M$. In this case we have
$$
\l_\lambda^{n+1}(1)(z_{n+1})=\l_\lambda(\l_\lambda^n(1))(z_{n+1})\le
\l_\lambda^n(1)(z_n)\l_\lambda(1)(z_{n+1})\le \l_\lambda^n(z_n)\cdot
$$
Those inequalities are implied by our assumptions on $M$ and $r$ that
ensure us first, that any pre-image of $z_{n+1}$ is in $\b_r$, and
second,
that $\l_\lambda(1)(z_{n+1})\le 1$. We may now apply our inductive
assumption to conclude that $H_{n+1}$ is true so that this hypothesis
is true for any integer $n$. Let $\alpha_{r,M,\lambda_0}$ be defined
as the infimum of the set $(m_\lambda(\b_{r,M}))$ where
$\lambda\in[1,\lambda_0]$. Since $\lambda\mapsto m_\lambda(\b_{r,M})$
is continuous on $[1,\lambda_0]$, this infimum is achieved
and is strictly greater than $0$. Fix $r$ small and choose $M(r)$ such
that all assumptions are fulfilled and set
$C_{r,\lambda_0}=\frac{K_{r,M(r)}}{\alpha_{r,M(r),\lambda_0}}$. We
deduce from our analysis that
$\lim_{n\to\infty}\l_\lambda^n(1)=\rho_\lambda\le C_{r,\lambda_0}$ on $\b_r$. We
have thus proved the left hand side inequality of Lemma
\ref{ineqconfinv}.
In order to finish the proof of this Lemma we need to prove that
$\frac{1}{K}\le\rho_\lambda\le K$ on $\b_{r,M}$. By (\ref{ineqconfinv01})
this will be done if one can prove that $\mu_\lambda(\b_{r,M})\ge \frac{1}{2}$ for suitable $r$ and $M$.

Since we know that $\rho_\lambda\le C_{r,\lambda_0}$ on $\b_r$, we may
already use the left-hand side inequalities of
Proposition~\ref{invariantmeasure}.
In particular for any $n$ we have
$$\mu_\lambda(C_n)\le
\frac{C_{r,\lambda_0}}{n^{2\delta-1}}\qquad\mbox{with
  $1<\delta=\inf\{d(\lambda)\}$, well defined by continuity.}$$
Let now $N$ be big enough so that
$$\sum_{n\ge N}\frac{1}{n^{2\delta-1}}\le \frac{1}{4C_{r,\lambda_0}}\cdot$$
Chose $r'$ small enough so that for any $\lambda\in[1,\lambda_0]$ we have
$$B(0,r')\subset \cup_{n\ge N} C_n(\lambda)\cdot$$
Such a choice is possible because of Proposition \ref{unifestimatesC}. We then easily conclude that $\mu_\lambda(B(0,r'))\le \frac{1}{4}$.
As a consequence, one may assume, without loss of generality, that we have started our analysis with $r>0$ small enough so that
$\mu_\lambda(B(0,r))\le \frac{1}{4}$.

By Lemma 4.1 in \cite{uz1}, we know that the sequence of measures $(m_\lambda)$ is tight. In particular, if $M$ is chosen large enough, then for
any $\lambda\in[1,\lambda_0]$ we have $m_\lambda(P_M^c)\le \frac{1}{4C_{r,\lambda_0}}$. From where we deduce that
$\mu_\lambda(P_M^c)\le\frac{1}{4}$.

Note now that $\mu_\lambda(\b_{r,M})=1-\mu_\lambda(B(0,r))-\mu_\lambda(P_M^c)\ge \frac{1}{2}$. As already mentioned this inequality finishes the
proof of the Lemma.
\cqfd
\section{Controlling the integrals}
In this section we mainly reproduce the reasoning of
\cite{hz}. Nevertheless there are some differences we would
like to emphasize : the main being that we do not know
whether the dimension of $J(F_1)$ is less than $\frac{3}{2}$ or not.
Note also that the Markov partition used in \cite{hz} is replaced in
the present article by
the backward images of the fundamental domain $C_0$. Finally, note that
we work directly on $J_\lambda$ without conjugating the dynamics.

Before we start the proofs and
in order to simplify some expressions and calculations, we
introduce the following notation. Let
$$
\Psi_n=\sum_{k=1}^{n}\frac{1}{(F_\lambda^k)'},
$$
$$
\Phi_n=\sum_{k=1}^{n}\frac{1}{|(F_\lambda^k)'|},
$$
$$
\Psi=\sum_{k=1}^{\infty}\frac{1}{(F_\lambda^k)'},
$$
and
$$
\Phi=\sum_{k=1}^{+\infty}\frac{1}{|(F_\lambda^k)'|}
$$
so that formula (\ref{formula}) may be written
$$
d'(\lambda)=-\frac{d(\lambda)}{\chi_{\mu_\lambda}}\left(1-\frac{1}{\lambda}\right)\int_{J_\lambda}
\re\left(\Psi\right)d\mu_\lambda\cdot
$$
We will need the following equation which is an easy computation
\begin{equation}\label{PsiPhi}
\Psi=\frac{1}{(F_\lambda^n)'}\Psi\circ F_\lambda^n+
\Psi_n,\qquad \Phi=\frac{1}{|(F_\lambda^n)'|}\Phi\circ F_\lambda^n+ \Phi_n\cdot
\end{equation}
\subsection{Lyapunov exponents}
In this paragraph we prove that the Lyapunov exponents do not play any
role in our estimates of the derivative. In order  to do
this we only need to check that they are uniformly bounded above and
separated away from zero. More precisely we prove the following.
\begin{Prop}\label{Lyapunov}
There exist $r_0>0$, $\lambda_0>1$ and $K>1$ such that $\forall \lambda\in (\lambda,\lambda_0)$ we have
$$\frac{1}{K}\le \chi_{\mu_\lambda}:=\int_{J_\lambda}\log |F_\lambda'|d\mu_\lambda\le K\cdot$$
\end{Prop}
\proof
First note that $\forall \lambda\ge 1$ and $\forall z\in J_\lambda$ we have $|F_\lambda'(z)|\ge 1$. In particular we have
$$\int_{C_0}\log |F_\lambda '|d\mu_\lambda \le \chi_{\mu_\lambda}\cdot$$
There is $K_1>0$ such that $\re z\ge K_1$ for any $z\in C_0(\lambda)$ and any $\lambda\in (1,\lambda_0)$, and by Proposition \ref{unifestimatesfan}
there is $K_2$ such that $\mu_\lambda(C_0)\ge K_2$. Since $\log |F_\lambda'(z)|=\log \lambda+\re z$ we deduce that
$$0<K_1K_2\le\int_{C_0}\log |F_\lambda '|d\mu_\lambda\le \chi_{\mu_\lambda}\cdot$$
This is the first part of the proof.
\vskip0.2cm
\noindent
For the other part note first that continuity of $\lambda\mapsto
d(\lambda)$ and the fact that $d(1)>1$ imply that there exist
$\alpha >1$ and $\beta>0$
such that $\alpha+\beta\le d(\lambda)$ for any $\lambda\in (1,\lambda_0)$. This implies in particular that $\forall \lambda\in (1,\lambda_0)$ and $\forall z\in J_\lambda$
$$\frac{1}{|(F_\lambda)'(z)|^{d(\lambda)}}\le\frac{1}{|(F_\lambda)'(z)|^{\alpha+\beta}}\cdot$$
Consider now the following partition of the strip $P$ : $(A_n)_{n\in
  \N}$, with $A_n:=\{z\in P\,|\, n-1<\re z\le n\}$. We have
$$
\begin{aligned}
\chi_{\mu_\lambda}
&=\sum_{n=1}^{+\infty}\int_{A_n}\log |F_\lambda'|d\mu_\lambda
\le \log \lambda_0+
\sum_{n=1}^{+\infty}\int_{A_n}\re z d\mu_\lambda(z) \\
&\le \log
\lambda_0+\sum_{n=1}^{+\infty}n\mu_\lambda(A_n)\cdot
\end{aligned}
$$
Lemma~\ref{ineqconfinv} implies that there exists $K_3>0$ such that
$\mu_\lambda(A_n)\le K_3m_\lambda(A_n)$ for $n\ge 2$. Note now that
$$m_\lambda (A_n)=\int_{J_\lambda}\chi_{A_n}dm_\lambda=\int_{J_\lambda}\l_\lambda(\chi_{A_n})dm_\lambda\cdot$$
For any $z\in J_\lambda$ and any $k\in \Z$ we let $z_k$ be the
preimage of $z$ for $F_\lambda$ such that $f_{\lambda}
(z_k)=z+2ik\pi$. We thus have
$$\l_\lambda(\chi_{A_n})(z)=\sum_{k\in\Z}\frac{1}{|F_\lambda'(z_k)|^{d(\lambda)}}\chi_{A_n}(z_k)\cdot$$
With $\alpha$ and $\beta$ defined above, this gives that
$$\l_\lambda(\chi_{A_n})(z)\le \sum_{k\in\Z}\frac{1}{|F_\lambda'(z_k)|^{\alpha+\beta}}\chi_{A_n}(z_k)\cdot$$
Since  $|F_\lambda'(z_k)|=\lambda e^{\re z_k}=|z+\lambda+2ik\pi|$, we have
$$\frac{1}{|F_\lambda'(z_k)|^{\alpha+\beta}}\chi_{A_n}(z_k)\le
\frac{1}{|z+\lambda+2ik\pi|^\alpha}\lambda^{-\beta}e^{-\beta (n-1)},
$$
so that
$$\l_\lambda(\chi_{A_n})(z)\le \lambda^{-\beta}e^{-\beta n}\sum_{k\in\Z}\frac{1}{|z+\lambda+2ik\pi|^\alpha}\cdot$$
As we have $\alpha>1$, there is $K_4>0$, independent of $\lambda$ and $z$, such that
$$\sum_{k\in\Z}\frac{1}{|z+\lambda+2ik\pi|^\alpha}<K_4\lambda^\beta\cdot$$
This tells us that
$$\l_\lambda(\chi_{A_n})(z)\le K_4e^{-\beta n}\cdot$$
Integrating with respect to $m_\lambda$, and summing over $n\ge 2$, we get
$$\chi_{\mu_\lambda}\le \log\lambda_0+K_3m_\lambda(A_1)+K_3K_4\sum_{n\ge 2}e^{-\beta n}\le K_5,$$
With $K_5:=\log\lambda_0+K_3+K_3K_4\frac{e^{-2\beta}}{1-e^{-\beta}}$. This is clearly independent of $\lambda$ and we are done
\cqfd

\

Note that with some more work one can indeed prove that
$\chi_{\mu_\lambda}$ converges towards $\chi_{\mu_1}$ as $\lambda$
converges towards $1$ from above.

\subsection{Controlling the integral away from $0$}
Let $N$ be an integer\footnote{This integer will be chosen later big enough to ensure that for any $z_n\in C_n$ we have
$\sum_{n\ge N}\arg(z_n)\le\frac{\pi}{6}$.} and set $M_N=\bigcup_{n\ge N+1} C_n$ and $B_N=J_\lambda\setminus M_N$. Note that both set $M_N$ and
$B_N$ depends on $\lambda$.
\begin{Prop}\label{intBN}
There exists $k(N)>0$ such that $\forall \lambda\in[1,\lambda_0]$ we have
\\$\ds \int_{B_N}\Phi d\mu_\lambda\le k(N)$.
\end{Prop}
\proof
Let $D_0=B_N$ and for any $n\in\N$ let $D_n=C_{N+n}$.
Following \cite{hz} let $U_n$ be the set of points which arrive or come back to $B_N$ after exactly $n$ iterates, which means
that $U_n=F_\lambda^{-1}(D_{n-1})$. Note that $U_n\cap M_n=D_n$. Given $N_0\in \N$ we set $A_n=F_\lambda^{-N_0}(U_n)\cap B_N$.
Since $(U_n)$ is a partition of $J_\lambda$, $(A_n)$ is a partition of $B_N$ and we have
\\$\ds \int_{B_N}\Phi d\mu_\lambda=\sum_{k=1}^{+\infty}\int_{A_k}\Phi d\mu_\lambda.$
\\Using relation \ref{PsiPhi} with $n=N_0+k$ we get
\\$\ds\int_{A_k}\Phi d\mu_\lambda=\ds\int_{A_k}\left(\frac{1}{|(F_\lambda^{N_0+k})'|}\Phi\circ F_\lambda^{N_0+k}+ \Phi_{N_0+k}\right)d\mu_\lambda\cdot$
\\Using the fact that $F_\lambda^{N_0+k}(A_k)\subset B_N$, Lemma \ref{landingoutsideB} and Lemma \ref{l2062805}  we deduce that
\\$\ds\int_{A_k}\Phi d\mu_\lambda\le \frac{\kappa(N)}{(N_0+k)^2}\ds\int_{A_k}\Phi\circ F_\lambda^{N_0+k}d\mu_\lambda+(N_0+k)\mu_\lambda(A_k)$
\\The fact that $F_\lambda^{N_0+k}(A_k)\subset B_N$ also implies that $\chi_{A_k}\le \chi_{B_N}\circ F_\lambda^{N_0+k}$, from
the invariance of $\mu_\lambda$ we thus get
\\$\ds \int_{A_k}\Phi\circ F_\lambda^{N_0+k}d\mu_\lambda\le \int\chi_{B_N}\circ F_\lambda^{N_0+k}\Phi\circ F_\lambda^{N_0+k}d\mu_\lambda
\le \int_{B_N}\Phi d\mu_\lambda$ which leads to
\\$\ds\int_{A_k}\Phi d\mu_\lambda\le \frac{\kappa(N)}{(N_0+k)^2}\ds\int_{B_N}\Phi d\mu_\lambda+(N_0+k)\mu_\lambda(A_k).$
\\In order to estimate $\mu_\lambda(A_k)$, we first use Lemma \ref{ineqconfinv} to conclude that $\mu_\lambda(A_k)\le K
m_\lambda(A_k)\le K m_\lambda (F_\lambda^{-N_0}{(U_k)})$, for some
constant $K$ independent of $k$, $N_0$ and $\lambda$. Since $U_k=F_\lambda^{-1}(D_k)$, we get
$\mu_\lambda(A_k)\le K m_\lambda (F_\lambda^{-(N_0+1)}{(D_k)})$.
Moreover
\\$\ds m_\lambda({F_\lambda^{-(N_0+1)}{(D_k)}})=\int \chi_{D_k}\circ F_\lambda^{N_0+1}dm_\lambda=\int_{D_k}\l_\lambda^{N_0+1}(1)dm_\lambda,$
\\since there exists $K_1(N_0)$ independent of $\lambda$ and $k$ such that $\ds\l_\lambda^{N_0+1}(1)\le K_1(N_0)$, using Lemma \ref{unifestimatesconf}
and the fact that $D_k=C_{N+k}$, we get
\\$\ds m_\lambda({F_\lambda^{-(N_0+1)}{(D_k)}}\le K_1(N_0)m_\lambda(C_{N+k})\le \frac{K_2}{(N+k)^{2d(\lambda)}}\cdot$
\\Using the fact that $(N_0+k)\le N_0(N+k)$, we thus conclude that
\\$\ds\int_{A_k}\Phi d\mu_\lambda\le \frac{\kappa(N)}{(N_0+k)^2}\ds\int_{B_N}\Phi d\mu_\lambda+\frac{KK_2N_0}{(N+k)^{2d(\lambda)-1}}\cdot$
\\Summing over $k$ we end up with
\\$\ds \int_{B_N}\Phi d\mu_\lambda\le \frac{\kappa(N)}{N_0-1}\ds\int_{B_N}\Phi d\mu_\lambda+\frac{KK_2N_0}{(N-1)^{2d(\lambda)-2}}
\cdot$
\\The integer $N$ being fixed, one may now choose $N_0$ big enough so that $\frac{\kappa(N)}{N_0-1}\le frac{1}{2}$, so that
\\$\ds \int_{B_N}\Phi d\mu_\lambda\le \frac{2KK_2N_0}{(N-1)^{2d(\lambda)-2}}
\cdot$
This last constant depends only on $N$ and we are done.
\cqfd
\subsection{Controlling the integral in a neighborhood of $0$}
In this paragraph we deal with the remaining part of $\int \re(\Psi) d\mu_\lambda$. If we  note $M_N=J_\lambda\setminus B_N$ we prove
\begin{Prop}\label{intMN}
There exists $K>0$ and $N\in \N$ such that for $\forall \lambda\in(1,\lambda_0)$
$$\frac{1}{K} (\lambda-1)^{2d(\lambda)-3}\le\int_{M_N}\re(\Psi)d\mu_\lambda\le K (\lambda-1)^{2d(\lambda)-3},\quad \mbox{if $d(\lambda)< \frac{3}{2}$,}$$
$$-\frac{1}{K} \log (\lambda-1)\le \int_{M_N}\re(\Psi)d\mu_\lambda\le -K \log(\lambda-1),\quad \mbox{if $d(\lambda)= \frac{3}{2}$,}$$
$$\left|\int_{M_N}\re(\Psi)d\mu_\lambda\right|\le K,\quad\mbox{if $d(\lambda)> \frac{3}{2}$.}$$
\end{Prop}
\proof
We split this integral into several pieces. First we note using \ref{PsiPhi} that
$$
\int_{M_N}\re(\Psi)d\mu_\lambda=
\sum_{n=N+1}^{+\infty}\left[\int_{C_n}\re\left(\frac{1}{(F_\lambda^{n-N})'}\Psi\circ F_\lambda^{n-N}\right)d\mu_\lambda
+\int_{C_n}\re\left(\Psi_{n-N}\right)d\mu_\lambda\right].
$$
We first deal with the left hand side of the sum that we bound integrating the modulus of the function.
$$\left|\int_{C_n}\re\left(\frac{1}{(F_\lambda^{n-N})'}\Psi\circ F_\lambda^{n-N}\right)d\mu_\lambda\right|\le
\int_{C_n}\frac{1}{|(F_\lambda^{n-N})'|}\Phi\circ F_\lambda^{n-N}d\mu_\lambda\cdot$$
We use Lemma \ref{landingoutsideB} and the fact that for $z\in C_n$, we have $F_\lambda^{n-N}(z)\in C_N\subset B_N$ to conclude that
$$\left|\int_{C_n}\re\left(\frac{1}{(F_\lambda^{n-N})'}\Psi\circ F_\lambda^{n-N}\right)d\mu_\lambda\right|\le \frac{K}{(n-N)^2}\int_{B_N}\Phi d\mu_\lambda\cdot$$
Summing over $n\ge N$ we get
$$\left|\sum_{n=N+1}^{+\infty}\left[\int_{C_n}\re\left(\frac{1}{(F_\lambda^{n-N})'}\Psi\circ F_\lambda^{n-N}\right)d\mu_\lambda\right]\right|\le
K\int_{B_N}\Phi d\mu_\lambda\sum_{n=1}^{+\infty}\frac{1}{n^2}\cdot$$
By Proposition \ref{intBN} we conclude that there exists $K(N)>0$ such that
\begin{eqnarray}\label{intMN2}
\left|\sum_{n=N}^{+\infty}\left[\int_{C_n}\re\left(\frac{1}{(F_\lambda^{n-N})'}\Psi\circ F_\lambda^{n-N}\right)d\mu_\lambda\right]\right|\le K(N)\cdot
\end{eqnarray}
We now deal with the right hand side. We have
$$\int_{C_n}\re(\Psi_{n-N})d\mu_\lambda=\sum_{k=1}^{n-N}\int_{C_n}\re\left(\frac{1}{(F_\lambda^k)'}\right)d\mu_\lambda\cdot$$
Choose $N$ big enough so that conclusions of Lemma \ref{Argumentznp} hold. For any $z\in C_n$ and any $k\le n-N$ we have
$$\frac{\sqrt{3}}{2}|F_\lambda^k(z)|\le\re (F_\lambda^k)'(z),$$
so that
$$\int_{C_n}\re(\psi_{n-N})d\mu_\lambda\sim \sum_{k=1}^{n-N}\int_{C_n}\re\left(\frac{1}{|F_\lambda^k|'}\right)d\mu_\lambda\cdot$$
Note now that for any $z\in C_n$, we have by the Chain Rule that
$$(F_\lambda^k)'(z)=\frac{(F_\lambda^n)'(z)}{(F_\lambda^{n-k})'(F_\lambda^k(z))},$$
with $F_\lambda^k(z)\in C_{n-k}$. We deduce, using Proposition \ref{unifestimatesCp}, that
$$\frac{1}{|(F_\lambda^k)'(z)|}\sim \frac{a_n(2)}{a_{n-k}(2)}\cdot$$
Estimates of $\mu_\lambda (C_n)$ are given by Proposition \ref{invariantmeasure} and we conclude that
$$\int_{C_n}\re(\psi_{n-N})d\mu_\lambda\sim\left\{\begin{array}{cc}
a_n(2d(\lambda)-1)a_n(2)\sum_{k=1}^{n-N} a_{n-k}(-2) & \mbox{if $n\le N_\varepsilon$,}
\vspace{0.2cm}\\
\frac{1}{\varepsilon}a_n(2d(\lambda))a_n(2)\sum_{k=1}^{n-N} a_{n-k}(-2)& \mbox{if $n\ge N_\varepsilon$.}
\end{array}\right.
$$
Since $a_n(\alpha)a_n(\beta)=a_n(\alpha+\beta)$, and with $S_{k,n}(\alpha)=\sum_k^n a_j(\alpha)$, this can also be written
$$\int_{C_n}\re(\psi_{n-N})d\mu_\lambda\sim\left\{\begin{array}{cc}
a_n(2d(\lambda)+1)S_{N,n-1}(-2) & \mbox{if $n\le N_\varepsilon$,}
\vspace{0.2cm}\\
\frac{1}{\varepsilon}a_n(2d(\lambda)+2)S_{N,n-1}(-2)& \mbox{if $n\ge N_\varepsilon$.}
\end{array}\right.
$$
Use now Corollary \ref{partialsum3} we have $S_{N,n-1}(-2)\sim (a_n(-3)-a_N(-3))$ if $n\le N_\varepsilon$ and $S_{N,n-1}(-2)\sim \frac{a_n(-2)}{\varepsilon}$
if $n\ge N_\varepsilon$ and we get
$$\int_{C_n}\re(\psi_{n-N})d\mu_\lambda\sim\left\{\begin{array}{cc}
a_n(2d(\lambda)+1)(a_n(-3)-a_N(-3)) & \mbox{if $n\le N_\varepsilon$,}
\vspace{0.2cm}\\
\frac{1}{\varepsilon^2}a_n(2d(\lambda)+2)a_n(-2)& \mbox{if $n\ge N_\varepsilon$.}
\end{array}\right.
$$
Since $a_n(\alpha)a_n(\beta)=a_n(\alpha+\beta)$ we get
$$
\int_{C_n}\re(\psi_{n-N})d\mu_\lambda\sim\left\{\begin{array}{cc}
a_n(2d(\lambda)-2)-a_N(-3)a_n(2d(\lambda)+1) & \mbox{if $n\le N_\varepsilon$,}
\vspace{0.2cm}\\
\frac{1}{\varepsilon^2}a_n(2d(\lambda)) & \mbox{if $n\ge N_\varepsilon$.}
\end{array}\right.
$$
Summing over $n\ge N$ this gives us $\sum_{n\ge N}\int_{C_n}\re(\psi_{n-N})d\mu_\lambda$ is comparable with
$$
\max\left((S_{N,N_\varepsilon}(2d(\lambda)-2)-a_N(-3)S_{N,N_\varepsilon}(2d(\lambda)+1)), \frac{1}{\varepsilon^2}S_{N_\varepsilon,+\infty}(2d(\lambda))\right)
$$
We then deduce from Corollary \ref{partialsum2} and Corollary \ref{partialsum3} that $S_{N,N_\varepsilon}(2d(\lambda)+1)\sim a_{N}(2d(\lambda))\sim 1$,
and also that $S_{N_\varepsilon,+\infty}(2d(\lambda))\sim \frac{a_{N\varepsilon}(2d(\lambda))}{\varepsilon}\sim \varepsilon^{2d(\lambda)-1}$. Estimates of
$S_{N,N_\varepsilon}(2d(\lambda)-2)$ depend on the comparison of $d(\lambda)$ with $\frac{3}{2}$. More precisely, if $d(\lambda)>\frac{3}{2}$ then
Corollary \ref{partialsum3} tells us that $S_{N,N_\varepsilon}(2d(\lambda)-2)\sim 1$, if $d(\lambda)=\frac{3}{2}$ then it tells us that
$S_{N,N_\varepsilon}(2d(\lambda)-2)\sim \log N_\varepsilon$, and if $d(\lambda)<\frac{3}{2}$ then $S_{N,N_\varepsilon}(2d(\lambda)-2)\sim \varepsilon^{2d(\lambda)-3}$.
Summarizing all those estimates we get
$$\sum_{n\ge N}\int_{C_n}\re(\psi_{n-N})d\mu_\lambda\sim
\left\{\begin{array}{cl}
1 & \mbox{if $d(\lambda)>\frac{3}{2}$}\\
\log N_\varepsilon  &  \mbox{if $d(\lambda)=\frac{3}{2}$} \\
\varepsilon^{2d(\lambda)-3} & \mbox{if $d(\lambda)<\frac{3}{2}$}
\end{array}
\right.$$

\subsection{Proof of the main result}
We are now in position to prove the main result of this paper that we
recall here.
\begin{Th}\label{main}
There exists $\lambda_0>1$, and $K>1$ such that
$$\left\{
\begin{array}{lclclc}
\frac{-1}{K}(\lambda-1)^{2d(1)-2} & \le  & d'(\lambda) & \le & -K(\lambda-1)^{2d(1)-2} & \mbox{if $d(1)<\frac{3}{2}$,}\\\hspace{0.1cm}
& & |d'(\lambda)| & \le & K (\lambda-1)\log \frac{1}{\lambda-1} & \mbox{if $d(1)=\frac{3}{2}$,}\\\hspace{0.1cm}
& & |d'(\lambda)| & \le & K(\lambda-1) & \mbox{if $d(1)>\frac{3}{2}$.}
\end{array}
\right.$$
In particular the function $\lambda\mapsto d(\lambda)$ is C$^1$ on $[1,+\infty)$, with $d'(1)=0$.
\end{Th}
\proof
Let us recall that we have
$$d'(\lambda)=-\frac{d(\lambda)}{\chi_{\mu_\lambda}}\left(1-\frac{1}{\lambda}\right)\int_{J_\lambda}
 \re\Psi d\mu_\lambda\cdot$$

We first use \cite{uz1}, where it is proved that $\lambda\mapsto d(\lambda)$ is continuous on $[1,+\infty)$
, and Proposition \ref{Lyapunov} to conclude that there exists $\lambda_1>1$ and $K_1>1$ such that
$\forall \lambda\in (1,\lambda_1)$ we have
$$\frac{1}{K_1}(\lambda-1)\le \frac{d(\lambda)}{\chi_{\mu_\lambda}}\left(1-\frac{1}{\lambda}\right)\le K_1(\lambda-1)\cdot$$
Note that given any integer $N$ we have
$$\int_{J_\lambda}
 \re\Psi d\mu_\lambda=\int_{B_N}\re\Psi d\mu_\lambda+\int_{M_N}\re\Psi d\mu_\lambda,$$
so that
\begin{eqnarray}\label{lastest}
|d'(\lambda)|\le 2K_1(\lambda-1)\max\left(\left|\int_{B_N}\re\Psi d\mu_\lambda\right|,\left|\int_{M_N}\re\Psi d\mu_\lambda\right|
\right)\cdot
\end{eqnarray}
We may thus use Proposition \ref{intBN} and Proposition \ref{intMN} to conclude that $d'(\lambda)$ is converging towards $0$
when $\lambda$ is converging towards $0$ from above. In particular there is $\lambda_2>1$ such that $\forall \lambda\in [1,\lambda_2)$,
$$-\frac{1}{2}\le d'(\lambda)\le \frac{1}{2}\cdot$$
We deduce that
$$-\frac{1}{2}(\lambda-1)\le d(\lambda)-d(1)\le \frac{1}{2}(\lambda-1),$$
so that
$$
\begin{aligned}
(\lambda-1)^{\lambda-1}(\lambda-1)^{2d(1)-3}
&\le(\lambda-1)^{2d(\lambda)-3}
=(\lambda-1)^{2d(1)-3}(\lambda-1)^{2(d(\lambda)-d(1))} \\
&\le (\lambda-1)^{-(\lambda-1)}(\lambda-1)^{2d(1)-3}
\end{aligned}
$$
Since $\lambda\mapsto (\lambda-1)^{\lambda-1}$ is continuous on $[1,\lambda_2]$ there exists $K_3>1$ such that
$$\frac{1}{K_3}(\lambda-1)^{2d(1)-3}\le (\lambda-1)^{2d(\lambda)-3}\le K_3(\lambda-1)^{2d(1)-3}\cdot$$
Using again Proposition \ref{intBN} and Proposition \ref{intMN}, and the fact we just proved that allows us to replace
$d(\lambda)$ with $d(1)$, we conclude the proof of the main result in case $d(1)<\frac{3}{2}$.
\vskip0.2cm
In case $d(1)=\frac{3}{2}$,  propositions \ref{intBN} and \ref{intMN} tells us that the maximum in (\ref{lastest})
is dominated by $-\log (\lambda-1)$. In case $d(1)>\frac{3}{2}$, the same proposition leads to the fact that this
maximum is bounded.
\cqfd

\section{Appendices}

\subsection{Estimates close to a repelling/parabolic fixed point}
In this appendix we show how to get estimates in
case  of a degeneracy towards a multi-petal
parabolic fixed point. It is a two steps
proof : first we deal with the real axis then we extend estimates
obtained in the real line to the complex plane using K\oe be's distortion
Theorem.

Consider the following family of germs of holomorphic functions
defined in a neighborhood of $0$ that we denote by $\mathcal{U}$:
$$
f_\varepsilon (z)=(1+\varepsilon)z+z^{p+1}+z^{p+2}g_\varepsilon
(z)\cdot
$$
Assume that there is an inverse branch $f_\varepsilon^{-1}$ well
defined on $\mathcal U$ that leaves a sector
$S_\theta:=\{re^{i\alpha}\,|\,-\theta\le\alpha\le\theta\}$ invariant,
for some $0<\theta<\frac{\pi}{2}$. Let $\u_\theta:=\u\cap S_\theta$. Assume also
that $\forall z\in\mathcal{U}$ we have $|zg_\varepsilon
(z)|<\frac{1}{2}$. Let $I=\mathcal{U}\cap\R^+$ and assume that
$f_\varepsilon^{-1} (I)\subset I$ and that $f_\varepsilon$ is not decreasing on
$I$.

This appendix is organized as follow : in the first two paragraphs we
study those germs giving in the second paragraph uniform estimates for
$|(f_{\varepsilon}^{-n})'|$.
\subsubsection{The mean value Theorem and its consequences}
We start with the following easy fact.
\begin{Lemme}\label{tech1}
Let $f:\R\to\R^+$ be a decreasing map with antiderivative $F$ on $\R$ and let
$(u_n)_{n\in\N}$ be a decreasing sequence of real numbers. Suppose that
there exist $n>1$ such that for all $k\le n$ we have $$
\begin{array}{lclc} \mbox{i-} & K_1\le (u_k-u_{k+1})f(u_k), &
\mbox{then} & K_1k\le F(u_0)-F(u_k)\cdot\\ \mbox{ii-} &
(u_k-u_{k+1})f(u_{k+1})\le K_2,&\mbox{then}& F(u_0)-F(u_k)\le K_2
k\cdot \end{array} $$
\end{Lemme}
\proof One only needs to check that our
assumptions imply $$(u_k-u_{k+1})f(u_k)\le\int_{u_{k+1}}^{u_k}f(t)dt\le(u_k-u_{k+1})f(u_{k+1})\cdot$$
\cqfd
\vskip0.2cm
In particular we point out the following two particular cases :
\begin{Cor}\label{cor1}
Let $(x_n)$ be a decreasing sequence of positive real numbers. Assume that
there exist $0<K_1<K_2$ and $n\in\N$ such that $\forall k\le n$,
$$K_1x_n^{p+1}\le (x_n-x_{n+1})\le K_2x_{n+1}^{p+1}.$$
Then there exist $\tilde{K_1}$ and $\tilde{K_2}$\footnote{One can take
 for instance $\tilde{K_2}=(pK_1)^{-\frac{1}{p}}$ and
$\tilde{K_1}=(pK_2+\frac{1}{x_0^p})^{-\frac{1}{p}}$.} such that for
$\forall k\le n$
$$\tilde{K_1}\le k^\frac{1}{p}x_k\le \tilde{K_2}\cdot
$$
\end{Cor}
\begin{Cor}\label{cor2}
Let $(u_n)$ be a decreasing sequence of real numbers. Assume that
there are $\alpha>0$, $\beta>0$, $p>0$  and $n\in\N$ such that $\forall k\le n$
$$(u_k-u_{k+1})\le \alpha+\beta e^{pu_{k+1}}\cdot$$
Then $\forall k\le n$ we have
$$\frac{\alpha^{\frac{1}{p}}}{(\alpha+\beta e^{pu_0})^\frac{1}{p}}e^{-\alpha k}\le e^{u_k-u_0}\cdot$$
\end{Cor}
Let us provide a short argument of how these corollaries can be
deduced from the Lemma~\ref{tech1}.
\proof
For Corollary \ref{cor1} we use the Lemma with the function $f:x\mapsto
x^{-(p+1)}$ so that one may take $F:x\mapsto -\frac{1}{p}x^{-p}$.
We deduce that we have :
$$K_1n\le\frac{1}{p}\left(\frac{1}{x_n^p}-\frac{1}{x_0^p}\right)\le
K_2n\cdot$$
Elementary computations then lead to the desired inequalities.
\vskip0.2cm
For Corollary \ref{cor2} we now consider the function
$f:x\mapsto (1+\frac{\beta}{\alpha}e^{px})^{-1}$. One first checks that
$F:x\mapsto x-\frac{1}{p}\log f(x)$ is an antiderivative of $f$.
Our assumptions on $(u_n)$ may now be written as
$$(u_k-u_{k+1})f(u_k)\le \alpha\cdot$$
Using the Lemma \ref{tech1} we deduce that $F(u_0)-F(u_k)\le \alpha k$.
This can be written in the form
$$
u_0-u_k+\frac{1}{p}\log\left(\frac{\alpha+\beta
    e^{pu_k}}{\alpha+\beta e^{pu_0}}\right)\le\alpha k\cdot$$
Applying exponents to both sides of this last inequality, we deduce that
$$\left(\frac{\alpha+\beta e^{pu_k}}{\alpha+\beta e^{pu_0}}\right)^{\frac{1}{p}}e^{-\alpha
k}\le e^{u_k-u_0}\cdot$$
From this we get our estimates.
\cqfd
\vskip0.2cm
\subsubsection{Uniform estimates along the real axis}
We now come back to our dynamical setting. Let $x_0\in I$ be a fixed
element. Assume for convenience that $x_0<1$.
Define for any $n\ge 0$, $f_\varepsilon(x_{n+1}(\varepsilon))=x_n({\varepsilon})$,
where $x_0(\varepsilon)=x_0$. For each $\varepsilon>0$ sufficiently small,
we define $N_\varepsilon$ as $N_\varepsilon=\sup\{n\in\N\,|\,x_n^p\ge
\varepsilon \},$ and for $\varepsilon=0$ as $N_0=+\infty$.
Note that for any $\varepsilon>0$ small enough, the sequence
$(x_n(\varepsilon))$ is strictly decreasing towards 0. So that
$N_\varepsilon$ is a well defined  integer. Our main results in this
paragraph is the following.
\begin{Prop}\label{unifestimates} There exists $K>1$ such that for all
 $\varepsilon>0$ small enough,
\begin{equation}\label{estN}
K^{-1}\le \varepsilon N_\varepsilon\le K,\quad \varepsilon>0,
\end{equation}
\begin{equation}\label{estx2}
K^{-1}\le x_nn^{\frac{1}{p}}\le K,\quad
\forall n< N_{\varepsilon},
\end{equation}
\begin{equation}\label{estx3}
K^{-1}\le x_n\varepsilon^{-\frac{1}{p}}(1+\varepsilon)^{n}\le K,\quad \forall
n\ge N_\varepsilon.
\end{equation}
\end{Prop}
This result may be interpreted in the following way : $N_\varepsilon$
is a "parabolic time". During that time, the fixed point $0$ acts
on the orbit of $x_0$, $(x_n)$, as if it was a parabolic fixed point
with $p$ petals. For $n$ greater than $N_\varepsilon$ the orbit
of $x_0$ is close enough to $0$ and realize that it is indeed an
attracting fixed point for $f_\varepsilon^{-1}$.

In the following Lemma we obtain estimates which are true for all $n\in\N$ and part of proposition \ref{unifestimates}.
\begin{Lemma}\label{unifestimatesfan} There exists $K>1$ such that for all
 $\varepsilon>0$ small enough,
\begin{equation}\label{estx1}
K^{-1}\varepsilon^{\frac{1}{p}}(1+\varepsilon)^{-n}\le x_n\le K
\frac{1}{n^\frac{1}{p}},\quad \forall n\in \N\cdot
\end{equation}
\end{Lemma}

\proof All our estimates will result from the following very definition of $(x_n)$.
\begin{equation}\label{defxn}
x_n=(1+\varepsilon)x_{n+1}+x_{n+1}^{p+1}(1+x_{n+1}g_\varepsilon(x_{n+1})).
\end{equation}
Assuming that $\varepsilon<1$, we easily deduce from this equality that for any $n$ we have
\begin{equation}\label{A1}
1\le\frac{x_n}{x_{n+1}}\le (2+2x_0^p)\le 4
\end{equation}
From \ref{defxn} we deduce that
\begin{equation}\label{A2}
\frac{x_n-x_{n+1}}{x_{n}^{p+1}}=\frac{\varepsilon x_{n+1}}{x_{n}^{p+1}}+\left(\frac{x_{n+1}}{x_n}\right)^{p+1}(1+x_{n+1}g_\varepsilon(x_{n+1})),
\end{equation}
which, with (\ref{A1}) leads to
$$K_1:=\frac{1}{4^{p+2}}\le\frac{1}{2}\left(\frac{x_{n+1}}{x_n}\right)^{p+1}\le\frac{x_n-x_{n+1}}{x_{n}^{p+1}}.$$
Using now Corollary \ref{cor1} we conclude that $\forall n$,
$x_n n^{\frac{1}{p}}\le \left(p\frac{1}{4^{p+2}}\right)^{\frac{-1}{p}}\le 64$.
This is precisely the right hand side (\ref{estx1}).

The left hand side of \ref{estx1} is obtained when one notes that
(\ref{defxn}) also implies that
$$\log x_n-\log x_{n+1}= \log
(1+\varepsilon+x_{n+1}^p+x_{n+1}^{p+1}g_{\varepsilon}(x_{n+1}))\le \log(1+\varepsilon) +2x_{n+1}^{p}.$$
We may thus apply Corollary \ref{cor2} with the sequence $u_n:=\log
x_n$,  $\alpha=\log(1+\varepsilon)$, and $\beta=2$, and deduce that
$$\frac{\alpha^{\frac{1}{p}}}{(\alpha+\beta e^{pu_0})^\frac{1}{p}}e^{-\alpha n}\le
e^{u_n-u_0}\cdot$$
Assuming that $\varepsilon$ is small enough so that $\frac{\varepsilon}{3}\le\alpha=\log(1+\varepsilon)\le e^{pu_0}$ we get
$$\frac{1}{9}\varepsilon^{\frac{1}{p}}(1+\varepsilon)^{-n}\le e^{u_n}=x_n.$$
This ends the proof of lemma \ref{unifestimatesfan}.
\cqfd
We are now in position to give a proof of Proposition \ref{unifestimates}, but first note that the right hand side of (\ref{estx2}) and
the left hand side of (\ref{estx3}) are given by Lemma \ref{unifestimatesfan}.
\proof
In order to get estimate (\ref{estx2}), we check that the assumptions on
$g_\varepsilon$, the definition of $N_\varepsilon$ and relation (\ref{A2})
leads for all $n<N_\varepsilon$ to
$$\frac{x_n-x_{n+1}}{x_{n+1}^{p+1}}\le \frac{5}{2}\cdot$$
Corollary \ref{cor1} then tells us that $\forall n<N_\varepsilon$
we have
$$\tilde{K_1}\le x_n n^{\frac{1}{p}},$$
with, for instance, $\tilde{K_1}=(\frac{5p}{2}+\frac{1}{x_0^p})^{-\frac{1}{p}}.$

We are now in position to give estimates for $N_\varepsilon$. They
easily come out from the following inequalities we have already
proved:
\begin{equation}
\frac{\tilde{K_1}}{K_0^{2p}}\frac{1}{N_\varepsilon^{\frac{1}{p}}}\le
\frac{\tilde{K_1}}{K_0^{2p}}\frac{1}{(N_\varepsilon-1)^{\frac{1}{p}}}\le
\frac{x_{N_\varepsilon-1}}{K_0^{2p}}\le
\frac{x_{N_\varepsilon}}{K_0}\le x_{N_\varepsilon+1}\le\varepsilon^{\frac{1}{p}}\le
x_{N_\varepsilon}\le \frac{2}{N_\varepsilon}\cdot
\end{equation}
From there we deduce that
\begin{equation}\label{A3}
\tilde{K_3}\le \varepsilon^{\frac{1}{p}}N_\varepsilon\le 2,
\end{equation}
with $\tilde{K_3}=\frac{\tilde{K_1}}{K_0^{2p}}$.

Now we only need to take care of (\ref{estx3}). We start by noticing that
for all $n$ we have $(1+\varepsilon)x_{n+1}\le x_n$. For any $n\ge
N_\varepsilon$ we thus have $(1+\varepsilon)^{n-N_\varepsilon}x_n\le
x_{N_\varepsilon}$. This leads to
\begin{equation}\label{A4}
x_n\le
(1+\varepsilon)^{-n}x_{N_\varepsilon}(1+\varepsilon)^{N_\varepsilon}.
\end{equation}
By definition of $N_\varepsilon$ and relation (\ref{A1}) we have
$$x_{N_\varepsilon}K_0 x_{N_\varepsilon+1}\le K_0\varepsilon^{\frac{1}{p}}\cdot$$
By relation (\ref{A3}) we also have
$$
(1+\varepsilon)^{N_\varepsilon}\le
(1+\frac{2^p}{N_\varepsilon})^{N_\varepsilon}
\le e^{2^p}\cdot
$$
From this and (\ref{A4}) we deduce that
\begin{eqnarray}
x_n\le K_0e^{2^p}\varepsilon^{\frac{1}{p}}(1+\varepsilon)^{-n}\cdot
\end{eqnarray}
Taking $K=K_0^{2p}e^{2^p}$ finishes the proof of the Proposition.

\cqfd

Let now $\alpha(p):=\frac{p+1}{p}$. The following corollary is
useful
\begin{Cor}\label{unifestimates2}
There exists $K>1$ such that for all $\varepsilon>0$ small enough
we have
$$
\begin{array}{rrclr}
\mbox{i-} &  K^{-1}  \le &
  (x_n-x_{n+1})n^{\alpha(p)}  & \le K & \forall n\le N_\varepsilon \\
\mbox{ii-} & K^{-1}  \le &
  (x_n-x_{n+1})\varepsilon^{-\alpha(p)}(1+\varepsilon)^{(p+1)n}
& \le K & \forall n
> N_\varepsilon \\
\end{array} $$
\end{Cor}
\proof
From relation (\ref{defxn}) we deduce that $\forall n\in\N$,
\begin{eqnarray}\label{unifestimates2A}
x_n-x_{n+1}=\varepsilon x_{n+1}+x_{n+1}^{p+1}(1+x_{n+1}g_\varepsilon(x_{n+1})),
\end{eqnarray}
so that $\varepsilon x_{n+1}\le x_n-x_{n+1}$, and Lemma~\ref{unifestimatesfan} tells us that the left hand side inequality of ii- is true.
Moreover, for $n\le N_\varepsilon$ we have $\varepsilon\le x_{n+1}^p$, and (\ref{unifestimates2A})
leads to $x_n-x_{n+1}\le 3x_{n+1}^{p+1}$. With Lemma~\ref{unifestimatesfan}, we get the right-hand side
inequality of i-.
\vskip0.2cm
Let $n<N_\varepsilon$. Then, from (\ref{estx2}) and (\ref{A1}), we deduce that
$$\frac{1}{2KK_0^{p+1}}\frac{1}{n^{\alpha(p)}}\le
\frac{x_n^{p+1}}{K_0^{p+1}}\le x_{n+1}^{p+1}\le (x_n-x_{n+1})\cdot$$
This is the left-hand side inequality of i-.
\vskip0.2cm
In order to finish fix $n\ge N_\varepsilon$. Then, by relation (\ref{defxn}) and
Proposition \ref{unifestimates}, we get that
$$x_n-x_{n+1}\le 3 x_{n+1}^{p+1}\le 3K^{p+1}\varepsilon^{\alpha(p)}
(1+\varepsilon)^{-(p+1)n}\cdot$$
The proof of the Corollary is now complete.
\cqfd
\subsubsection{Extension to the complex plane}
As already mentioned, this extension is done via K\oe be's
distortion Theorem. It asserts that given two simply
connected domains in $\C$, $V\subset V'$, such that the boundary
of $V$ is at a positive distance from the boundary of $V'$, there exists
a constant $K>0$, which depends only on the modulus of the
annulus $V'\setminus V$, and such that {\it for any univalent function
$f$ defined in $V'$} we have for all $x,y\in V$ we have,
$$
\frac{1}{K}
\le\frac{|f'(x)|}{|f'(y)|}\le K\cdot
$$

\begin{Prop}\label{unifestimatesC}
Let $V$ be a domain such that
$\bar{V}\subset\u_\theta$. Then there exists $K>0$ such
that $\forall \varepsilon$ small enough, $\forall n\in\N$ and
$\forall z\in V$ we have
$$
\begin{array}{rrllr}
\mbox{i-} & \frac{1}{K}  \le &
  n^{\alpha(p)}|(f_\varepsilon^{-n})'(z)| & \le K & \mbox{if
    $n<N_\varepsilon$.} \vspace{0.2cm} \\
\mbox{ii-} & \frac{1}{K}(1+\varepsilon)^{-(p+1)n}  \le &
 |(f_\varepsilon^{-n})'(z)|   & \le
  K\varepsilon^{\alpha(p)}(1+\varepsilon)^{-(p+1)n}& \mbox{if $n\ge N_\varepsilon$.} \\
\end{array} $$
\end{Prop}
\proof
Enlarging $V$ if necessary one may assume that there is $x_0\in
V\cap\R^+\cap\u_\theta$ such that for all
$\varepsilon$ small enough $x_1(\varepsilon):=f_\varepsilon^{-1}(x_0)$
is also in $V$. K\oe be's distortion Theorem implies that for all
$n$, all $\varepsilon$ and all $z\in V$ we have
$$\frac{1}{K}\frac{(x_n(\varepsilon)-x_{n+1}(\varepsilon))}{x_0-x_1(\varepsilon)}\le
|(f_\varepsilon^{-n})'(z)|\le
K\frac{(x_n(\varepsilon)-x_{n+1}(\varepsilon))}{x_0-x_1(\varepsilon)}\cdot$$
Applying Corollary \ref{unifestimates2}, and noticing that
$x_0-x_1(\varepsilon)>a>0$ with some real $a$ independent of $\varepsilon$,
lead to the desired inequalities.
\cqfd

\

The following result gives uniform estimates on how closely the orbits
are tangent to the real axis.

\begin{Cor}\label{Argumentzn}
There exists $K>0$ such  that $\forall \varepsilon$ small enough, $\forall n\in\N$ and
$\forall z_0\in V$,  we have
$$|\im (f_\varepsilon^{-n}(z_0))|\le K\frac{1}{n^{\alpha(p)}}\cdot$$
In particular, the series $\sum_{n=0}^\infty\im (f_\varepsilon^{-n}(z_0))$ converges.
\end{Cor}
\proof Note that $|\im(z_n)|=|\im(z_n-x_n)|\le |z_n-x_n|$. K\oe be's
distortion theorem leads to $|z_n-x_n|\le
K\frac{1}{|(f_\varepsilon^{n})'(z_n)|}$ and  Proposition
\ref{unifestimatesC} gives the result.
\cqfd

\subsection{Estimates of some partial sums}
In this appendix we single out the behaviour of the partial sums we
need to evaluate at several steps in the proof of our main result.
It seemed to us that postponing those estimates to an appendix will
clarify the exposition. We are thus in this paragraph dealing
with a sequence of real numbers defined by : $a_n=1/n$ for
$n\le N_\varepsilon$ and $a_n=\varepsilon (1+\varepsilon)^{-n}$ for
$n>N_\varepsilon$, where $N_\varepsilon$ is comparable with
$1/\varepsilon$. We are indeed interested in the sequences
$(a_n(\alpha))_{n\in\N}$, with $\alpha\in\R$ and
$a_n(\alpha)=a_n^\alpha$, and partial sums
$S_{k,n}(\alpha)=\sum_{j=k}^{n}a_j(\alpha)$.

The first Lemma, whose proof is straightforward and left to the reader
asserts, the following.

\begin{Lemma}\label{partialsum1}
For any $k<n$ in $\N$ we have
$$S_{k,n}(\alpha)\sim \left\{
\begin{array}{cl}
\frac{1}{1-\alpha}(n^{1-\alpha}-k^{1-\alpha}) & \mbox{if $n\le N_\varepsilon$ and $\alpha\ne 1$}\\
\log\frac{n}{k} & \mbox{if $n\le N_\varepsilon$ and $\alpha= 1$}
\end{array}\right.
$$
$$
S_{k,n}(\alpha)\sim \frac{1}{\alpha\varepsilon}(a_k
(\alpha)-a_n(\alpha)),\,\mbox{if $k>N_\varepsilon$
and $\alpha\ne 0$.}$$
\end{Lemma}
As its consequence, we get the following.
\begin{Cor}\label{partialsum2}
If $\alpha>0$ then
$$\begin{array}{ccccl}
\mbox{i-} & S_{n,+\infty}(\alpha) & \sim & \frac{a_n(\alpha)}{\varepsilon} & \mbox{if $n>N_\varepsilon$,}
\\
\mbox{ii-} & S_{n,+\infty}(\alpha) & \sim & a_n(\alpha-1) & \mbox{if $n\le N_\varepsilon$ and $\alpha> 1$,}
\\
\mbox{iii-} & S_{n,+\infty}(\alpha) & \sim &
\log\frac{N_\varepsilon}{n}+K & \mbox{if $n\le N_\varepsilon$,
  $\alpha=1$, for some $K>0$.}
\\
\mbox{iv-} & S_{n,+\infty}(\alpha) & \sim & N_\varepsilon^{1-\alpha}
& \mbox{if $n\le N_\varepsilon$ and $\alpha< 1$,}
\end{array}
$$
\end{Cor}
\proof Since $\alpha >0$, we see that the sequence $(1+\varepsilon)^{-\alpha n}$ converges
to $0$, and Lemma~\ref{partialsum1} implies that  i- is true. Note that we have
$$\max(S_{n,N_\varepsilon}(\alpha),S_{N_\varepsilon,+\infty}(\alpha))\le S_{n,+\infty}\le
2\max(S_{n,N_\varepsilon}(\alpha),S_{N_\varepsilon,+\infty}(\alpha)).$$
Using i- that we have just proved, the fact that  we have
$a_{N_\varepsilon}\sim a_{N_\varepsilon+1}$, and the fact that
$N_\varepsilon\sim\varepsilon^{-1}$,
we conclude that
$$S_{N_\varepsilon,+\infty}\sim
\frac{a_{N_\varepsilon}}{\varepsilon}\sim \varepsilon^{\alpha-1}\sim
N_\varepsilon^{1-\alpha}\cdot$$
Let us now estimate $S_{n,N_\varepsilon}$ by considering three cases.
We start with the case when $\alpha=1$. Indeed, Lemma \ref{partialsum1}
implies that $S_{n,N_\varepsilon}\sim \log(\frac{N_\varepsilon}{n})$. This
gives us  iii-.

Assume now that $\alpha>1$. Then $S_{N_\varepsilon,+\infty}\sim
N_\varepsilon^{1-\alpha}\le n^{1-\alpha}=a_n(\alpha-1)$. Moreover,
in virtue of Lemma \ref{partialsum1},
we have $S_{n,N_\varepsilon}\sim n^{1-\alpha}-N_\varepsilon^{1-\alpha}$. Thus
$$
S_{n,N_\varepsilon}\sim
a_n(\alpha-1)\left(1-\left(\frac{n}{N_\varepsilon}\right)^{\alpha-1}\right)\cdot
$$
In particular $S_{n,N_\varepsilon}\lesssim a_n(\alpha-1)$. So, we can
conclude that $S_{n,+\infty}(\alpha)\lesssim a_n(\alpha-1)$. If
$\frac{n}{N_\varepsilon}\le \frac{1}{2}$, we have
$(1-\frac{n}{N_\varepsilon})^{\alpha-1}
\ge (1-\frac{1}{2})^{\alpha-1}$. And we also have
$S_{n,+\infty}(\alpha)\gtrsim a_n(\alpha-1)$; so, we are done. On the
other hand, if $\frac{n}{N_\varepsilon}\ge \frac{1}{2}$, then
$$S_{n,+\infty}(\alpha)\ge S_{N_\varepsilon,+\infty}(\alpha)\sim
N_\varepsilon^{1-\alpha}\sim n^{1-\alpha}=a_n(\alpha-1)\cdot$$
This ends the proof of  ii-.

Assume finally that  $0<\alpha<1$. Then Lemma \ref{partialsum1} tells us that
$$
S_{n,N_\varepsilon}(\alpha)\sim
(N_\varepsilon^{1-\alpha}-n^{1-\alpha})\le
N_\varepsilon^{1-\alpha}\sim S_{N_\varepsilon,\infty}.$$
We thus conclude that
$\max(S_{n,N_\varepsilon}(\alpha),S_{N_\varepsilon,+\infty}(\alpha))\sim
S_{N_\varepsilon,+\infty}$. This proves iv-
and ends the proof of the Corollary.
\cqfd

\

We can also prove the following result with the same kind of
arguments. So we omit them.
\begin{Cor}\label{partialsum3}
Let $N$ be a fixed integer such that
$2N<N_\varepsilon\sim\frac{1}{\varepsilon}$.
Then we have the following estimates of $S_{N,n}(\alpha)$ for $N\le n$:
$$S_{N,n}(\alpha)\sim
\left\{
\begin{array}{ll}
\left.\begin{array}{cl}
a_N(\alpha-1)-a_n(\alpha-1) & \mbox{for $1<\alpha$} \\
\log\frac{n}{N} & \mbox{for $\alpha=1$}\\
a_n(\alpha-1)-a_N(\alpha-1)& \mbox{for $\alpha<1$}\\
\end{array}\right\} & \mbox{for $n\le N_\varepsilon$}\\
\vspace{0.2cm}
\left.\begin{array}{cl}
1 & \mbox{for $1<\alpha$}\\
\log N_\varepsilon & \mbox{for $1=\alpha$}\\
N_\varepsilon^{1-\alpha} & \mbox{for $0\le\alpha<1$}\\
\frac{a_n(\alpha)}{\varepsilon} & \mbox{for $\alpha<0$}\\
\end{array}\right\} & \mbox{for $n\ge N_\varepsilon$}
\end{array}\right.
$$
\end{Cor}

\end{document}